\documentclass[headings=standardclasses,
	11pt,
	DIV=12,
	bibliography=totoc
]{scrartcl}

\usepackage[T1]{fontenc}
\usepackage{authblk}
\usepackage{mathtools, amssymb, amsthm, mathrsfs, bm}
\usepackage{graphicx}
\graphicspath{{figures/}}
\usepackage[english]{babel}
\usepackage[shortlabels]{enumitem}
\usepackage{algorithm, algpseudocode, float}
\usepackage{tabularray}
\usepackage{microtype}
\usepackage{csquotes}
\usepackage[sort&compress]{natbib}
\usepackage{xcolor,ninecolors}
\usepackage[colorlinks=true, linkcolor=blue4, citecolor=green6, 
pdfstartview=FitH, pdfpagemode=UseNone]{hyperref}

\newcommand*{\F}{\mathcal F}

\newcommand*{\R}{\mathbb R}

\newcommand*{\V}{\mathcal V}

\newcommand*{\intR}{\int_\R}
\renewcommand*{\d}{\mathrm{d}}
\newcommand*{\sgne}{\mathrm{sign}_\epsilon}
\newcommand*{\dt}{\frac{\d}{\d t}}
\newcommand*{\pd}[1][t]{\frac{\partial}{\partial #1}}
\newcommand*{\pdd}[1][x]{\frac{\partial^2}{\partial {#1}^2}}

\newcommand*{\ek}{\bm{e}_k}
\newcommand*{\phv}{\, \cdot \,} % PlaceHolder Value
\newcommand*{\finf}{f_\infty}

\newcommand*{\x}{x_*}
\newcommand*{\xx}{X_*}
\newcommand*{\xa}{x_\gamma}
\newcommand*{\xxa}{X_\gamma}
\newcommand*{\X}{X_t^i}

\renewcommand*{\(}{\begin{equation}}
\renewcommand*{\)}{\end{equation}}

\DeclareMathOperator*{\argmin}{arg\,min}
\DeclarePairedDelimiter{\ev}{\langle}{\rangle}
\DeclarePairedDelimiter{\abs}{\lvert}{\rvert}

\DeclarePairedDelimiterXPP{\abse}[1]{}\lvert\rvert{_\epsilon}{#1}
\DeclarePairedDelimiterXPP{\normpq}[1]{}\lVvert\rVvert{_\alpha}{#1}

\theoremstyle{plain}

\theoremstyle{remark}
\newtheorem{remark}{Remark}

\title{Superlinear drift in consensus-based optimization with condensation phenomena}
\date{}
\author[1]{Jonathan Franceschi\thanks{\texttt{jonathan.franceschi@unife.it }}}
\author[2,1]{Lorenzo Pareschi\thanks{\texttt{L.Pareschi@hw.ac.uk }}}
\author[3]{Mattia Zanella\thanks{\texttt{mattia.zanella@unipv.it}}}

\affil[1]{\normalsize Department of Mathematics and Computer Science, \protect\\ University of Ferrara, Ferrara, Italy}
\affil[2]{\normalsize Maxwell Institute for Mathematical Sciences and Department of Mathematics,\protect\\ Heriot Watt University, Edinburgh, UK}
\affil[3]{\normalsize Department of Mathematics ``F. Casorati'', \protect\\ University of Pavia, Pavia, Italy}

\begin{document}
\maketitle

\begin{abstract}
Consensus-based optimization (CBO) is a class of metaheuristic algorithms designed for global optimization problems. In the many-particle limit, classical CBO dynamics can be rigorously connected to mean-field equations that ensure convergence toward global minimizers under suitable conditions.  
In this work, we draw inspiration from recent extensions of the Kaniadakis--Quarati model for indistinguishable bosons to develop a novel CBO method governed by a system of SDEs with superlinear drift and nonconstant diffusion. The resulting mean-field formulation in one dimension exhibits condensation-like phenomena, including finite-time blow-up and loss of $L^2$-regularity. To avoid the curse of dimensionality a marginal based formulation which permits to leverage the one-dimensional results to multiple dimensions is proposed.  
We support our approach with numerical experiments that highlight both its consistency and potential performance improvements compared to classical CBO methods.
\end{abstract}

\textbf{Keywords}: Consensus-based optimization, superlinear drift, Kaniadakis--Quarati model, mean-field limit, finite-time blow-up, condensation 

\tableofcontents

\section{Introduction}

In recent decades, interacting particle systems have emerged as powerful tools for investigating many-agent dynamics of complex systems characterized by self-organizing features. These models are particularly prominent in fields such as mathematical biology, ecology, and socio-economic sciences, where they have been applied to phenomena including swarming, crowd dynamics, opinion formation, and synchronization processes~\cite{albim3a2019,carrillo2010particle,HaTadmor,choi_etal,COUZIN20021,carrillo_toscani,MT,CPT,PTZ,CS,DOB,Barbaro}. 

In recent times,  these results have also inspired the development of algorithmic frameworks for solving high-dimensional and computationally demanding optimization problems. In particular, a variety of metaheuristic approaches, including evolutionary algorithms and swarm intelligence methods, derive their structure from the behavior of interacting agents. These algorithms are increasingly employed in global optimization tasks, such as the NP-hard traveling salesman problem, where traditional deterministic approaches often fall short in terms of scalability or efficiency. 

Metaheuristics are generally regarded as high-level strategies for exploring the parameter space. They aim to balance exploration and exploitation, often through stochastic mechanisms that allow for adaptive and robust performance across a range of problem classes. Among the most studied and applied methods are Genetic Algorithms (GA), Particle Swarm Optimization (PSO), Ant Colony Optimization (ACO), and Simulated Annealing (SA), see e.g.,~\cite{goldberg1989genetic,Dorigo2019,KE,KGV}. Despite their widespread use and  success in practical applications, many metaheuristics still lack rigorous mathematical foundations that would justify their performance guarantees or convergence properties. A central challenge in this area remains the theoretical validation of whether, and under what conditions, these algorithms can reliably approximate optimal solutions.

Along these lines, Consensus-Based Optimization (CBO) provides a robust mathematical foundation for the systematic study of metaheuristic approaches to global optimization problems of the form 
\[
\x \in \textrm{argmin}\, \mathcal F(x), \qquad x \in \mathbb R^d,
\]
where $\mathcal F(\phv)\colon \mathbb R^d \to \mathbb R$ is a given non-convex, high dimensional, cost function. The original CBO method was introduced in~\cite{pinnau2017consensus, carrillo2021consensus} and has been further analyzed in a variety of settings, including convergence and mean-field formulations~\cite{fornasier2024consensus,FHPS,huang2022meanfield,Carrillo2018AnalyticalCBO,huang2024uniform,ham3a2020,kom3a2022}.  
Since then, several extensions of the method have been developed.
In particular, variants for constrained optimization~\cite{Borghi2023ConstrainedCBO} and multi-objective optimization~\cite{Borghi2023AdaptiveMOCBO} have been proposed.  
We also note that the mean-field approach underlying the CBO framework has been generalized to second-order models, providing connections to PSO algorithms~\cite{grassi2021particle}. Similar arguments have lead to mean-field interpretations of GA and SA~\cite{AFT,Pareschi2024SimulatedAnnealing,borghi2025kinetic,HZ25}. For recent surveys on these topics, we refer the reader to~\cite{borghi2024kinetic,Totzeck2022}.

The core idea is to guide the evolution of realizations in the search space via a system of interacting stochastic differential equations. In particular, it is assumed that  $X^i_t \in \mathbb R^d$, $i = 1,\dots,N$ follows the trajectory
\(\label{eq:CBO}
\d X^i_t = - \lambda (X^i_t - X_\gamma)\d t + \sigma \sum_{k = 1}^d \ek(X^i_t - X_\gamma)_k\d B_k^i, 
\)
where, for any $k = 1,\dots,d$, the term $\{B_k^i\}_{i=1}^N$ is a family of independent Wiener processes, $\ek$ is the $k$th vector of the canonical basis of $\mathbb{R}^d$,  and 
\[
X_\gamma = \dfrac{\sum_{j=1}^N X^i_t e^{-\gamma \mathcal F(X^i_t)}}{\sum_{j=1}^N e^{-\gamma \mathcal F(X^i_t)}}.
\]

This choice is rooted on the Laplace principle such that for every compactly-supported probability measure $\psi \in \mathcal P(\R^d)$ and such that $\textrm{argmin} \F(x)$ belongs to the support of $\psi$, then
\[
\lim_{\gamma \to +\infty} \biggl(-\frac{1}{\gamma}
	\log\Bigl(\int_{\R^d} e^{-\gamma \F(x)}\, \d \psi(x) \Bigr)
	\biggr) = \min_{\operatorname{supp(\psi})} \F(x),
\]
so that when $\gamma \gg 1$ we can approximate $x_\gamma \approx \argmin\F(X)$. 

Formally, it has been shown that in the limit $N \to +\infty$ the empirical density $f_{N}(x,t) = \frac{1}{N} \sum_{i=1}^N \delta(x-X_i)$ associated to the particle system $\{X_i\}_{i=1}^N \in \R^d$ converges to the continuous density $f(x,t)$ solution to the following Fokker--Planck equation-type with nonconstant diffusion
\[
\pd f(x,t) = \lambda \nabla \cdot \bigl([x - \xa(t)] f(x,t)\bigr)
	+ \dfrac{\sigma^2}{2} \sum_{k = 1}^d \partial_{kk} \bigl([x - \xa(t)]^2_k f(x,t)\bigr),
\]
where
\(\label{eq:xadef}
\xa(t) \coloneqq
	\frac{\int_{\R^d} xe^{-\gamma\F(x)} f(x,t) \, \d x}
		{\int_{\R^d} e^{-\gamma\F(x)} f(x,t)\, \d x}.
\)
Furthermore, under suitable assumptions on the drift and diffusion parameters $\lambda$ and $\sigma$, the solution $f(x,t)\to \delta(x-x_*)$ for $t \to +\infty$ exponentially fast in time, being $x_*$ the global minimizer of $\mathcal F$. 

In contrast to the classical CBO method, where the alignment dynamics rely on a linear attraction toward the weighted consensus point ${x}_\gamma$, the present work introduces a new family of consensus-based optimization schemes characterized by a \emph{superlinear drift}. This modification is intended to enhance mass concentration near the global minimizer and accelerate convergence, while also introducing a rich phenomenology, including finite-time blow-up, that parallels condensation phenomena observed in nonlinear kinetic equations.
 The construction is inspired by recent developments in nonlinear Fokker--Planck models for indistinguishable bosons, and in particular by extensions of the Kaniadakis--Quarati framework
\cite{kaniadakis1993kinetic,toscani2012}. The resulting dynamics exhibit strong concentration effects and, under appropriate conditions, finite-time blow-up reminiscent of condensation phenomena. These features are particularly relevant in the context of global optimization, where sharp localization around minimizers is a desirable property.

We recall that the original work~\cite{kaniadakis1993kinetic} was proposed to describe the emergence of Bose--Einstein statistics, which in $\mathbb{R}^3$ exhibits a critical mass and finite-time blow-up in the supercritical regime. In this direction, the effect of superlinear drift terms has been investigated in~\cite{toscani2012,BAGT,CHR,carrillo20081d}, highlighting their role in the formation of singularities and condensate-like behaviour.
More recently, the interplay between superlinear drift and nonconstant diffusion in Fokker--Planck equations has been studied in the context of consensus-type dynamics~\cite{calzola2024emergence}, where condensation phenomena were observed. In the subsequent work~\cite{toscani2024condensation}, it was further shown that finite-time blow-up persists even in the presence of nonconstant diffusion, provided the initial datum is sufficiently regular. In particular, for any $f(x,0) \in L^1 \cap L^\infty$, we observe blow-up of the solution in finite time.

To illustrate the potential of the approach, we perform a detailed one-dimensional analysis. In this setting, the new model allows for a  study of the effects of the drift exponent $\alpha$ on the behaviour of the solution. Explicit results highlight the role of the critical mass in producing finite-time condensation around the global minimizer. The one-dimensional case also reveals the computational and analytical challenges induced by high values of $\alpha$, where finite-time blow-up and strong singularities emerge.
In multiple dimensions, however, the situation becomes more involved. The reconstruction of the particle density needed to evaluate the nonlinear drift becomes prohibitively expensive due to the curse of dimensionality. Moreover, a rigorous mean-field analysis in this setting remains out of reach, as the presence of anisotropic and concentrated interactions poses significant theoretical obstacles. To circumvent some of these issues, we propose a novel \emph{marginal-based extension} of the method. This reformulation preserves the key features of the one-dimensional dynamics while dramatically reducing the computational complexity. 

The remainder of the paper is organized as follows. In Section~\ref{sec:model}, we introduce the consensus-based optimization model featuring a superlinear drift, and establish its connection with a nonlinear Fokker--Planck equation with nonconstant diffusion. Section~\ref{sec:meanfield} is devoted to the one-dimensional case, focusing in particular on the emergence of condensation phenomena and finite-time blow-up. Next in Section~\ref{sec:marginal} we introduce a novel method based on the use of marginals which is suitable for multi-dimensional simulations. In Section~\ref{sec:numerics}, we present a comprehensive set of numerical experiments that both validate the theoretical predictions and benchmark the performance of the proposed method against the classical CBO approach. To this aim, in higher dimensions, we introduce the marginal-based extension of the algorithm. Finally, the last section provides concluding remarks and outlines perspectives for future research.

\section{A superlinear CBO method}\label{sec:model}

Following the notations introduced in~\cite{pinnau2017consensus} and~\cite{carrillo2021consensus}, we consider a particle system $\X \in \mathbb R^d$, $j = 1,\ldots, N$, whose evolution is governed by the following system of SDEs
\(\label{eq:general-sde}
\d \X = -\lambda(\X - \xxa(t)) (1 + \beta K(\X) f_{\epsilon,N}^\alpha(\X))\d t
		+ \sigma \sum_{k = 1}^d \ek(H(\X))_k\d B_k^i,
\)
where $\ek$ is the $k$-th unit vector of the standard basis in~$\R^d$ and $B_k^i$ are $N\times d$ independent Wiener processes, while $\alpha,\beta \ge0$. In \eqref{eq:general-sde} we further introduced two nonnegative functions $H,K\colon \mathbb R^d \to \mathbb R^+$. In more detail, $K(\cdot)$ takes into account a state-dependent drift towards $\xxa$, while $H(\cdot)$ governs the ampliture of the noise, possibly allowing for anisotropic diffusion.  We may observe that \eqref{eq:general-sde} corresponds to the classical CBO under the hypothesis $\beta = 0$ and $H(x) = (x-X_\gamma)_k$. In \eqref{eq:general-sde} we introduced the mollifier $f_{\epsilon,N}$ of the atomic measure $f_N$, which is defined as
\begin{equation}
\label{eq:fN_reg}
f_{\epsilon,N}(x,t) \coloneqq \dfrac{1}{N} \sum_{j=1}^N \psi_\epsilon(x-\X).
\end{equation}
The function $\psi_\epsilon(\cdot)$ is a mollifier. Classical choices proposed in the literature are given by the Gaussian mollifier, such that
\[
\psi_\epsilon(y) \coloneqq \dfrac{1}{(2\pi\epsilon)^{d/2}}\exp\left\{ -\dfrac{\abs{y}^2}{2\epsilon}\right\}, 
\]
for any $\epsilon>0$, see~\cite{CARRILLO2020100066} and the references therein, or the piecewise-constant mollifier, defined as  
\[
\psi_\epsilon(y) = \dfrac{1}{\epsilon^d}\chi(y \in B_\epsilon),
\]
where $\chi(\cdot)$ denotes the indicator function and $B_\epsilon$ is the $d-$dimensional hypercube having bin width $\epsilon>0$. Proceeding as illustrated in~\cite{CP}, we have that $ f_{\epsilon,N}(x,t)  \overset{\ast}{\rightharpoonup} f_N(x,t) \in \mathcal P(\mathbb R^d)$, being $f_N(x,t)$ the empirical distribution of the  particle system.

To illustrate the features of the superlinear drift we fix $\xxa= a\in \mathbb R^d$ to be a constant vector. Therefore, we get
\begin{equation}
\label{eq:power2}
\begin{split}
\dfrac{\d}{\d t} \mathbb E[\abs{\X-a}^2] &= -2\lambda \mathbb E[\abs{\X-a}^2] - 2\lambda \beta \mathbb E[\abs{\X - \xxa}^2K(\X)f_{\epsilon,N}^\alpha]  + \sigma^2 \sum_{i=1}^d \mathbb E[\abs{\X-a}^2]_i \\
&\le -2\lambda \mathbb E[\abs{\X-a}^2]+ \sigma^2 \sum_{i=1}^d \mathbb E[\abs{\X-a}^2]_i 
\end{split}
\end{equation}
Therefore, the decay of the shifted second order moment is obtained again under the condition
\begin{equation}
\label{eq:assumption}
\sigma^2  < 2\lambda,  
\end{equation}
consistently with what obtained in \cite{carrillo2021consensus}. 
However, since $K(\cdot)\ge0$ and $f_{\epsilon,N}\ge0$, equation~\eqref{eq:power2} shows that, for any $\beta>0$, the density--dependent term induces an additional dissipative contribution. As a consequence, the condition~\eqref{eq:assumption} is sufficient but no longer sharp, and consensus formation may occur under milder assumptions compared to the regime $\beta=0$. At the numerical level, this allows to consider a wider range of values for the diffusion coefficient $\sigma^2$, thereby enabling a broader exploration of the parameter space. 

Indeed, the nonnegativity of both $K(\cdot)$ and $f_{\epsilon,N}$ results in an effective amplification of the drift towards the consensus point $\xxa(t)$. Consequently, particles located in regions of higher empirical density experience a stronger attraction, which promotes consensus with respect to the case $\beta=0$.

\subsection{Mean-field limit}
When the number of agents become large it is convenient to describe the evolution of the particle system as a mean-field equation. We limit ourselves to few classical general reference on this topics 
\cite{Golse,CIP} and several approaches to consensus and flocking hydrodynamics~\cite{Toscani,pinnau2017consensus,carrillo_toscani,carrillo2010particle}. 

We consider the regularization of the atomic probability measure $f_{\epsilon,N}(x,t)$ defined in \eqref{eq:fN_reg} and associated to the solution of the system \eqref{eq:general-sde}.  
For any $i = 1,\dots,N$ and for any test function $\varphi \in \mathcal C_c^\infty(\mathbb R^d)$, from the multidimensional It\^o's formula we get
\begin{equation}
\label{eq:ito_law}
\begin{split}
d\varphi(\X)
	&= \nabla\varphi(\X) \cdot \left( -\lambda(\X - \xxa)(1+\beta K(\X)f_{\epsilon,N}^\alpha(\X) \right)\d t \\
	&\hphantom{{}=}+ \sigma \sum_{k = 1}^d \ek H(\X) \nabla\varphi(\X)\ dB_k^i +\frac{\sigma^2}{2}  \sum_{k = 1}^d \ek H^2(\X) \nabla^2 \varphi(\X):I\, \d t,
\end{split}
\end{equation}
where $\nabla^2\varphi$ is the Hessian and $A:B = \textrm{Tr}(A^TB)$. Therefore, summing all the contributions of the particle system, since 
\[
\begin{split}
\dfrac{1}{N} \sum_{i=1}^N \d\varphi(\X)
	&= 	\dfrac{1}{N}\sum_{i=1}^N\nabla\varphi(\X) \cdot \left( -\lambda(\X - \xxa)(1+\beta K(\X)f_{\epsilon,N}^\alpha(\X) \right)\d t \\
	&\hphantom{{}=}+ \dfrac{\sigma}{N} \sum_{i=1}^N\ \sum_{k = 1}^d \ek H(\X) \nabla\varphi(\X)\d B_k^i\\
	&\hphantom{{}=}+ \dfrac{\sigma^2}{2N} \sum_{i=1}^N\sum_{k = 1}^d \ek H^2(\X) \nabla^2 \varphi(\X):I\, \d t,
\end{split}
\]
which gives
\[
\partial_t \ev*{f_N,\varphi} = \ev*{f_N, \nabla\varphi(x)\cdot\lambda(x-x_\gamma)(1+\beta K(x)f_{\epsilon,N}^\alpha(x))}  + \ev[\bigg]{f_N, \frac{\sigma^2}{2} \sum_{k=1}^d \ek H^2(x)\nabla^2\varphi(x):I} 
\]
Hence, after integration by parts, we obtain
\[
\ev[\bigg]{\partial_t f_N - \nabla\cdot \bigl[\lambda(x-x_\gamma) f_N (1+\beta K(x) f^\alpha_{\epsilon,N}(x))\bigr] - \frac{\sigma^2}{2}\sum_{k=1}^d \partial_{kk}\bigl(H_k(x)^2 f_N\bigr) ,\varphi}= 0
 \]
 or, in strong form 
 \begin{equation}
 \label{eq:strong}
 \partial_t f_N =  \lambda\nabla\cdot \bigl[ (x-x_\gamma) f_N (1+\beta K(x) f^\alpha_{\epsilon,N}(x))\bigr] +\frac{\sigma^2}{2}\sum_{k=1}^d \partial_{kk}(H_k(x)^2 f_N).
 \end{equation}
For each $t \ge 0$, the variance of the particle system \eqref{eq:general-sde} is decreasing in time under the assumption~\eqref{eq:assumption}. As a consequence, the empirical density $f_N$ is a probability measure in $\mathcal P(\mathbb R^d)$ with support uniformly bounded with respect to $N$. By Prokhorov's theorem,  the sequence $\{f_N\}_{N\in \mathbb N}$ is weakly-$\ast$-relatively compact, and there exists a subsequence $\{f_{N_\ell}\}_{\ell>0}$ and a mesure $f \in \mathcal P(\mathbb R^d)$ such that
\[
f_{N_\ell} \overset{\ast}{\rightharpoonup} f, \qquad \ell\to +\infty. 
\]
From \eqref{eq:strong}, in the limit $\ell\to +\infty$ and $\epsilon\to 0$ of the subsequence $\psi_\epsilon * f_{N_\ell}$, we formally get

\(\label{eq:FP-multiD}
\dfrac{\partial}{\partial t}f(x, t) = \lambda\nabla \cdot \Bigl[(x - \xa(t)) f(x, t)
						\bigl(1 + \beta K(x) f^\alpha(x,t)\bigr)  \Bigr] 
						+ \frac{\sigma^2}{2} \sum_{k = 1}^d \dfrac{\partial}{\partial x_{kk}} \bigl(H^2_k(x) f(x,t)\bigr),
\)
where we indicate with $(\phv)_k$ the $k$-th component of $(\phv)$.

In kinetic models for bosons, such as the Kaniadakis--Quarati or
Bose--Einstein Fokker--Planck equations, it is known that the formation
of finite-time singularities depends on a critical balance between the
 exponent $\alpha\ge0$ and the spatial dimension $d \ge 1$.
In particular, condensation phenomena may occur when
$\alpha > \frac{2}{d}$, whereas for $\alpha \le \frac{2}{d}$ the
diffusion is sufficiently strong to prevent the formation of
singularities; see, for instance, \cite{CHR,BAGT}.

Although proving an analogous result for our superlinear consensus model \eqref{eq:FP-multiD} is challenging due to its nonlocal structure, this observation suggests the existence of a dimension-dependent critical threshold which may have a negative impact on the application of the model to high dimensional optimization problems. This motivates the one-dimensional analysis of Section \ref{sec:meanfield} and the subsequent marginal-based formulation discussed in Section~\ref{sec:marginal}.

\section{Finite time blow-up in one-dimension}\label{sec:meanfield}

Throughout the one-dimensional analysis for $d = 1$, we will rely on the following assumptions:

\begin{enumerate}[(H1)]
  \item\label{hypothesis:1} The function $K(x)$ is nonnegative, and $H^2(x)\in \mathcal C^2(\mathbb{R})$.
  \item\label{hypothesis:2} The initial empirical measure $f_{N}(x,0)$ converges to a probability density $f_0\in L^1(\mathbb{R})\cap L^\infty(\mathbb{R})$.
  \item\label{hypothesis:3} The mollified empirical density $f_{\epsilon,N}(x)$ in \eqref{eq:fN_reg} converges in $L^1_{\mathrm{loc}}$ to a limiting density $f \in \mathcal{C}([0,T]; \mathcal{P}(\mathbb{R}))$, where $\mathcal{P}(\mathbb{R})$ denotes the set of probability measures on $\mathbb{R}$.
  \item\label{hypothesis:4} The global consensus point $x_\gamma(t)$ is defined as
  \[
  x_\gamma(t) = \frac{\int_{\mathbb{R}} x\, e^{-\gamma F(x)} f(x,t)\, dx}{\int_{\mathbb{R}} e^{-\gamma F(x)} f(x,t)\, dx}.
  \]
\end{enumerate}

In the case $d= 1$ we can easily see that equation~\eqref{eq:FP-multiD}  reads as follows
\(\label{eq:FP1Dmov}
	\pd f(x,t) = \lambda\pd[x] \bigl[[x - \xa(t)]f(x,t)(1 + \beta K(x) f^\alpha(x,t)\bigr]
	+ \frac{\sigma^2}{2} \pdd \bigl[H^2(x)f(x,t)\bigr].
\)

From \eqref{eq:FP1Dmov} we get, for any test function $\varphi(x)$,
\[
\begin{aligned}
\pd \int_{\mathbb R}\varphi(x) f(x,t)\, \d x
	&=  - \lambda\int_{\mathbb R} \frac{\d}{\d x} \varphi(x) \Bigl[(x-x_\gamma(t))f(x,t)(1+\beta K(x) f^\alpha(x,t))\\
	&\hphantom{{}=}+ \frac{\sigma^2}{2} \pd[x] \bigl(H^2 f(x,t)\bigr) \Bigr]\,\d x,
\end{aligned}
\]
and the mass is conserved whereas we cannot conclude the preservation of higher order moments due to the nonlinearity in the drift term. 

Our current requirements on the coefficients $K(x)$ and $H(x)$ are sufficient to ensure the positivity of the solution to~\eqref{eq:FP1Dmov}. Indeed, we can leverage a classical technique employed e.g., in~\cite{carrillo20081d,escobedo1991large}. We show that the $L^1(\R)$ of the solution $f(x,t)$ to~\eqref{eq:FP1Dmov} is non-increasing. This, along with mass conservation, immediately implies that the positive part of $f(x,t)$ coincides with $f(x,t)$ itself at all times. The proof entails taking an increasing regularized approximation of the sign function, $\sgne(\phv)$, so that its first derivative $\sgne'(\phv) > 0$. Moreover, we call $\abse{\phv}$ the primitive of $\sgne(\phv)$. This way, if we multiply equation~\eqref{eq:FP1Dmov} by $\sgne(f)$ and integrate by parts with respect to $x$ (omitting the dependence of $f$ on $x$ and $t$ for brevity), we get
\[
\begin{aligned}
\dt \intR \abse{f}\, \d x
	&= -\lambda\intR \sgne' f \pd[x] f (x - \xa(t))(1 + \beta K(x)f^\alpha)\, \d x\\
	&\hphantom{{}=}	-\frac{\sigma^2}{2} \intR \sgne' f \pd[x] f \Bigl(\pd[x]H^2(x) f + H^2(x)\pd[x]f\Bigr)\, \d x\\
	&= -\lambda\intR \pd[x]\Bigl[f \sgne f -\abse{f}\Bigr]
		\Bigl((x - \xa(t))(1 + \beta K(x)f^\alpha + \frac{\sigma^2}{2} \pd[x]H^2(x)\Bigr)\,\d x\\
	&\hphantom{{}=} -\frac{\sigma^2}{2}\intR \sgne' f \Bigl(\pd[x]f\Bigr)^2 H^2(x)\, \d x,
\end{aligned}
\]
where the second equality comes from observing that
\[
\pd[x]\bigl(f\sgne f - \abse{f} \bigr) = f\sgne' f \pd[x]f,
\]
so that
\[
\begin{aligned}
\dt \intR \abse{f}\, \d x &= \lambda\intR \Bigl[f \sgne f -\abse{f}\Bigr]
	\pd[x]\Bigl((x - \xa(t))(1 + \beta K(x)f^\alpha + \frac{\sigma^2}{2}\pd[x]H^2(x)\Bigr)\,\d x\\
	&\hphantom{{}=} -\frac{\sigma^2}{2}\intR \sgne' f \Bigl(\pd[x]f\Bigr)^2 H^2(x)\, \d x,
\end{aligned},
\]
and taking the limit for $\epsilon \to 0^+$ we have
\(
\dt \intR \abs{f} \, \d x \le  0,
\)
proving the claim.

\subsection[\texorpdfstring{Loss of $L^2$-regularity and evolution of energy-type functionals}{Loss of L²-regularity and evolution of energy-type functionals}]{Loss of $\bm{L^2}$-regularity and evolution of energy-type functionals}

We can observe that, for any nonnegative function $f \in L^2(\mathbb R)$, the second order moment of $f$ about the point~$b\in\R$ controls the $L^1$ norm since for any $R>0$ we get the bound
\[
\begin{split}
\int_{\mathbb R} f(x)\, \d x
	&= \int_{\abs{x - b} \le R} f(x)\, \d x + \int_{\abs{x-b} > R} f(x)\, \d x\\
	&\le (2R)^{1/2} \biggl(\int_{\mathbb R} f^2(x)\, \d x \biggr)^2 + \dfrac{1}{R^2} \int_{\mathbb R} (x - b)^2 f(x)\, \d x
\end{split}
\]
and optimizing on $R\ge0$ we get that 
\[
\int_{\mathbb R} f(x) \, \d x \le 5 \biggl(\dfrac{1}{2\sqrt{2}} \biggr)^{4/5} \biggl(\int_{\mathbb R} f^2(x)\, \d x\biggr)^{2/5}  \biggl( \int_{\mathbb R} (x - b)^2f(x)\,  \d x\biggr)^{1/5}. 
\]
The above estimate yields a bound on the $L^1-$norm of a nonnegative function in terms of its $L^2-$norm and its second order moment. Therefore, the boundedness of the $L^2$ norm of $f$ prevents its second order moment from vanishing and in particular
\[
\int_{\mathbb R} (x - b)^2 f(x)\, \d x \ge \biggl( \dfrac{1}{5}\biggr)^5 (2\sqrt{2})^4 \dfrac{\displaystyle\biggl(\int_{\mathbb R} f(x)\, \d x\biggr)^5}{\displaystyle\biggl( \int_{\mathbb R} f^2(x)\, \d x\biggr)^2}. 
\]
Hence, if $f$ is the solution to \eqref{eq:FP1Dmov} we may exploit the relationship between its $L^2$ norm and the evolution of its second order moment about the point~$b$ to check the preservation of $L^2$~regularity.  Based on these observations, we follow~\cite{fornasier2024consensus} and we consider the evolution of 
\(\label{eq:Vdef}
\V(t) \coloneqq \intR (x - \x)^2 f(x,t)\, \d x, \qquad \x \coloneqq \argmin_{\R} \F(x).
\)
whose evolution is given by 
\(\label{eq:gendV1}
\begin{aligned}
\dt \V(t)
	&= -2\lambda\intR (x - \x)(x - \xa(t))f(x,t)\, \d x\\
	&\hphantom{{}=}- 2\lambda\beta \intR (x - \x)(x - \xa(t))K(x)f^{\alpha + 1}(x,t)\, \d x\\
	&\hphantom{{}=} + \sigma^2 \intR H^2(x) f(x,t)\, \d x.
\end{aligned}
\)
Now, we notice that
\[
(x - \x)(x - \xa(t)) = (x - \x)^2 - (x - \x)\bigl(\xa(t) - \x\bigr),
\]
which translates equation~\eqref{eq:gendV1} into
\(
\begin{aligned}
\dt \V(t)
	&= -2\lambda\V(t) + 2\lambda\bigl(\xa(t) - \x\bigr)\intR (x - \x) f(x,t)\, \d x 
	+\sigma^2 \intR H^2(x)f(x,t)\, \d x\\
	&\hphantom{{}=}- 2\lambda\beta \intR (x - \x)^2 K(x) f^{\alpha + 1}(x,t)\, \d x\\
	&\hphantom{{}=} + 2\lambda\beta \bigl(\xa(t) - \x\bigr)\intR (x - \x) K(x) f^{\alpha + 1} (x,t)\, \d x\\
\end{aligned}
\)
From Cauchy-Schwarz inequality (cfr.~\cite{fornasier2024consensus}), we have
\(
\begin{aligned}
\dt \V(t)
	&\le -2\lambda\V(t) + \sqrt{2}\lambda\abs{\xa(t) - \x}\sqrt{\rho}\sqrt{\V(t)}
	+\sigma^2 \intR H^2(x)f(x,t)\, \d x\\
	&\hphantom{{}=}- 2\lambda\beta \intR (x - \x)^2 K(x) f^{\alpha + 1}(x,t)\, \d x\\
	&\hphantom{{}=} + 2\lambda\beta \abs{\xa(t) - \x}\underbrace{\intR \abs{x - \x} K(x) f^{\alpha + 1} (x,t)\, \d x}_{(I)}.
\end{aligned}
\)
The main idea is that Laplace's principle should allow us to say that
\[
\xa(t) \xrightarrow[\gamma \to +\infty]{} \x.
\] 
Hence, provided  $(I) = o(|x_\gamma - \x|)$ we may simplify the obtained inequality as follows
\(\label{eq:simpledV}
	\dt \V(t)
		\le -2\lambda\V(t) + \sigma^2 \intR H^2(x)f(x,t)\, \d x
		- 2\lambda\beta \intR (x - \x)^2 K(x) f^{\alpha + 1}(x,t)\,\d x. 
\)

For analytical purposes, we consider the case $K(x)\equiv 1$ and $H(x) = \abs{x-\x}$.  With $K(x)$ constant, the nonlinear interaction term in \eqref{eq:FP-multiD} no longer depends explicitly on position,  allowing us to focus on the effect of the density-dependent term $f^\alpha$. 
We focus first on the rightmost term in \eqref{eq:simpledV}, i.e.,
\[
A(t) \coloneqq \intR (x - \x)^2 f^{\alpha + 1}(x,t)\, \d x.
\]
We introduce the notation $\Omega \coloneqq \{x \in \R :\abs{x - \x} < R\}$ for a given positive radius $R$, so that we can write
\(\label{eq:ptrick}
\begin{aligned}
\rho = \intR f(x,t)\, \d x
	&= \int_\Omega f(x,t)\, \d x + \int_{\Omega^C} f(x,t)\, \d x\\
	&\le \int_\Omega [(x - \x)^2]^{-p}[(x - \x)^2]^p f(x,t)\, \d x 
	+  \frac{1}{R^2} \int_{\Omega^C} (x - \x)^2 f(x, t)\, \d x,\\
	&\le \int_\Omega [(x - \x)^2]^{-p}[(x - \x)^2]^p f(x,t)\, \d x 
	+  \frac{1}{R^2} \intR (x - \x)^2 f(x, t)\, \d x,
\end{aligned}
\)
see~\cite{toscani2024condensation}.
The goal is to derive a suitable H\"older-type inequality. 
To guarantee integrability of $[(x - \x)^2]^{-k}$ on $\Omega$, we require $k < 1/2$.
Choosing $p = 1/(\alpha+1)$, we have
\[
\int_\Omega f(x,t)\, \d x = \int_\Omega [(x - \x)^2]^{-\frac{1}{\alpha+1}}[(x - \x)^2]^{\frac{1}{\alpha+ 1}} f(x,t)\, \d x,
\]
which becomes, choosing the H\"older conjugate pair $\alpha + 1$, $(\alpha + 1)/\alpha$
\[
\int_\Omega [(x - \x)^2]^{-\frac{1}{\alpha+1}}[(x - \x)^2]^{\frac{1}{\alpha+1}} f(x,t)\, \d x
\le (c_\alpha)^{\frac{\alpha}{\alpha+1}}R^{\frac{\alpha - 2}{\alpha + 1}}
\left(\int_\Omega (x - \x)^2 f^{\alpha + 1}(x,t)\, \d x\right)^{\frac{1}{\alpha + 1}},
\]
where the constant $c_\alpha>0$, for any $\alpha>2$, results from
\[
\int_\Omega [(x - \x)^2]^{-\frac{1}{\alpha}}\, \d x = \frac{2\alpha}{\alpha - 2}R^{\frac{\alpha - 2}{\alpha}} = c_\alpha R^{\frac{\alpha - 2}{\alpha}}<+\infty.
\]
Therefore, we get
\(\label{eq:g}
\rho \le
(c_\alpha)^{\frac{\alpha}{\alpha+1}}R^{\frac{\alpha - 2}{\alpha + 1}}
\left(\int_\Omega (x - \x)^2 f^{\alpha + 1}(x,t)\, \d x\right)^{\frac{1}{\alpha + 1}}
+ \frac{1}{R^2} \intR (x - \x)^2 f(x, t)\, \d x.
\)
Optimising with respect to the radius $R>0$ we obtain
\(\label{eq:Wbound2}
A(t) \ge \frac{\rho^{\frac{3\alpha}{2}}}{(c_\alpha d_\alpha + d_\alpha^{-2})^{\frac{3\alpha}{2}} \V^\frac{\alpha - 2}{2}(t)},
\)
where, following~\cite{toscani2024condensation}, we have introduced the notation
\[
d_\alpha = [2(\alpha+1)/(c_\alpha(\alpha - 2))]^{1/3},\qquad \alpha>2. 
\]
Thanks to the bound introduced above, from \eqref{eq:simpledV} we obtain
\begin{equation}\label{eq:Vaftertrick}
\dt \V(t) \le (\sigma^2 - 2\lambda)\V(t) - 2\lambda\beta A(t)
	\le (\sigma^2 - 2\lambda)\V(t) - 2\lambda\beta A(0).
\end{equation}
Under the hypothesis $\sigma^2 \le 2\lambda$, as in~\cite{carrillo2021consensus}, the first term on the right-hand side is non-positive and can be disregarded. 
It then follows that
\begin{equation}
\V(t) \le \V(0)
	- 2\lambda\beta \frac{\rho^{\frac{3\alpha}{2}}}{(c_\alpha d_\alpha + d_\alpha^{-2})^{\frac{3\alpha}{2}} \V^{\frac{\alpha - 2}{2}}(0)} \, t.
\end{equation}
Therefore, $\bar t>0 $ exists such that  $\V(t) = 0$ for all $t \ge \bar t$. In particular, we get
\(\label{eq:finaltime}
\bar t = \frac{(c_\alpha d_\alpha + d_\alpha^{-2})^{\frac{3\alpha}{2}} \V^{\frac{\alpha}{2}}(0)}{2\lambda\beta\rho^{\frac{3\alpha}{2}}},
\)
which is positive for any $\alpha > 2$. Therefore, $f(x,t)$ solution to 
\[
\dfrac{\partial}{\partial t} f(x,t)  =  \dfrac{\partial}{\partial x} \left[ \lambda(x-x_\gamma)(1+\beta f^\alpha(x,t))f(x,t) + \frac{\sigma^2}{2} \dfrac{\partial}{\partial x} (\abs{x-\x}^2 f(x,t))\right], 
\] 
loses $L^2$ regularity in finite time for any $\rho>0$.

\subsection{Existence of a critical mass}
Let us generalize the choice in the previous subsection and take $H(x) = \abs{x - \x}^{\eta}$, with $0<\eta < 1/2$, and also $K(x) = H^{2\alpha}(x)$. If we fix $x_\gamma \equiv \x$, we can consider the following equation
\begin{equation}
\label{eq:sol2}
\dfrac{\partial}{\partial t} f(x,t)  =  \dfrac{\partial}{\partial x} \left[\lambda (x-\x)\bigl(1+\beta H^{2\alpha}(x)f^\alpha(x,t)\bigr)f(x,t) + \frac{\sigma^2}{2}\dfrac{\partial}{\partial x} \bigl(H^2(x) f(x,t)\bigr)\right], 
\end{equation}
which may be rewritten in terms of $g(x,t) = H^2(x)f(x,t)$ as follows
\[
\dfrac{1}{H^2(x)}\dfrac{\partial}{\partial t} g(x,t) = \dfrac{\partial}{\partial x} \left[ \lambda\dfrac{x-\x}{H^2(x)}g(x,t)(1+\beta g^\alpha(x,t))  + \frac{\sigma^2}{2}\dfrac{\partial}{\partial x} g(x,t)\right],
\]
whose unique steady state is the solution to 
\[
\dfrac{1}{1+\beta g^\alpha_\infty(x)}\dfrac{\partial}{\partial x}g_\infty(x)=  - \dfrac{4\lambda}{\sigma^2} \lambda\dfrac{x-\x}{H^2(x)}. 
\]
For any $\eta<1/2$ we get
\[
g_\infty^\alpha(x) = \dfrac{C^\alpha \exp\left\{ - \frac{2\lambda\alpha}{\sigma^2} \frac{(x-\x)^{2(1-\eta)}}{1-\eta} \right\}}{1-\beta C^\alpha\exp\left\{ - \frac{2\lambda\alpha}{\sigma^2} \frac{(x-\x)^{2(1-\eta)}}{1-\eta} \right\} }, 
\]
being $C>0$, and
\[
f_\infty(x) = \dfrac{1}{(x-\x)^{2\eta}} \dfrac{C \exp\left\{ - \frac{2\lambda}{\sigma^2} \frac{(x-\x)^{2(1-\eta)}}{1-\eta} \right\}}{\left(1-\beta C^\alpha\exp\left\{ - \frac{2\lambda\alpha}{\sigma^2} \frac{(x-\x)^{2(1-\eta)}}{1-\eta} \right\}\right)^{1/\alpha} }. 
\]
Consequently the constant $C>0$ has to satisfy the condition 
\[
C \exp\left\{ - \dfrac{2\lambda}{\sigma^2} \dfrac{(x-\x)^{2(1-\eta)}}{1-\eta}\right\}\le \left(\dfrac{1}{\beta}\right)^{1/\alpha}, 
\]
whose maximum is attained at $x = \x$ and implies that $C$ cannot be greater than
\[
C_M = \left(\dfrac{1}{\beta}\right)^{1/\alpha}. 
\]
Furthermore, we get
\[
\lim_{x \to \x} \dfrac{C \exp\left\{ - \dfrac{2\lambda}{\sigma^2} \dfrac{(x-\x)^{2(1-\eta)}}{1-\eta}\right\}}{(x-\x)^2} = \frac{\sigma^2}{2\lambda} \alpha, 
\]
and $f_\infty(x)$ diverges in $x = \x$ with order $2/\alpha$ as $x \to \x$. As a consequence, setting $f_{\infty,C_M}(x)$ the steady state obtained with the critical constant $C_M$,  the integral 
\[
\rho_c = \int_{\mathbb R} f_{\infty,C_M}(x)\, \d x < +\infty
\]
as soon as $\alpha>2$. Finally, the mass $\rho_c$ defines the critical mass of the system as for the classical Kaniadakis-Quarati model. 

Additionally, we can prove that the solution to \eqref{eq:sol2} loses regularity, in finite time, in the weighted space
\[
L^2_w = \Bigl\{f : \int_{\mathbb R} H^2(x) f^2 \, dx < +\infty \Bigr\}.
\] 
As a preliminary observation, we note that
\[
\intR H^2(x) f(x,t)\, \d x = \intR (x - \x)^{2\eta} f(x,t)\, \d x,
\]
and its evolution satisfies
\(\label{eq:Hevolution}
\begin{aligned}
\dt \intR H^2(x) f(x,t)\, \d x
	&= -2\eta\lambda \intR H^2(x) f(x,t)\, \d x - 2\eta\lambda\beta \intR H^{2\alpha+2}(x) f^{\alpha+1}(x, t)\, \d x\\
	&\hphantom{{}=} + \sigma^2\eta(2\eta - 1) \intR \Bigl(\pd[x] H^2(x) \Bigr)^2 f(x,t)\, \d x\\
	&\le  -2\eta\lambda \intR H^2(x) f(x,t)\, \d x - 2\eta\lambda\beta \biggl(\intR H^2(x) f(x, t)\, \d x\biggr)^{\alpha + 1}\\
	&\hphantom{{}=}+ \sigma^2\eta(2\eta - 1) \intR \Bigl(\pd[x] H^2(x) \Bigr)^2 f(x,t)\, \d x,
\end{aligned}
\)
where the inequality comes from Jensen inequality. 
The last term on the right-hand side is negative for $\eta < 1/2$, 
which ensures that the functional decays at least exponentially in time.
Similarly, for the related functional
\[
\intR H^4(x) f(x,t)\, \d x = \intR (x - \x)^{4\eta} f(x,t)\, \d x,
\]
its evolution is given by
\(\label{eq:H2evolution}
\begin{aligned}
\dt \intR H^4(x) f(x,t)\, \d x
	&= -4\eta\lambda \intR H^4(x) f(x,t)\, \d x - 4\eta\lambda\beta \intR H^{2\alpha+4}(x) f^{\alpha+1}(x, t)\, \d x\\
	&\hphantom{{}=} + 2\eta\sigma^2(4\eta - 1) \intR (x - \x)^{6\eta - 2} f(x,t)\, \d x,
\end{aligned}
\)
which decays in time if $\eta < 1/4$.

Equations~\eqref{eq:Hevolution} and~\eqref{eq:H2evolution} can be leveraged to analyze the evolution of another energy-type functional, i.e.,
\[
U(t) \coloneqq \intR (x - \x)^2 H^2(x) f(x,t)\, \d x = \intR (x - \x)^{2\eta + 2} f(x,t)\, \d x,
\]
whose evolution reads
\(\label{eq:Uevolution}
\begin{aligned}
	\dt U(t)
	&\le - \lambda(2\eta + 2) U(t) - (2\eta + 2)\lambda\beta \intR (x - \x)^2 H^{2\alpha + 2}(x) f^{\alpha + 1}(x, t)\, \d x\\
	&\hphantom{{}=} + \sigma^2(\eta + 1)(2\eta + 1) \intR H^4(x) f(x,t)\, \d x. 
\end{aligned}
\)
\begin{remark}
Notice that in the case $\eta = 0$, i.e., $H \equiv 1$, the functional $U(t)$ coincides with $\V(t)$ and equation~\eqref{eq:Uevolution} reads
\(\label{eq:Ueta0evolution}
\begin{aligned}
	\dt \V(t)
	&\le - 2\lambda \V(t) - 2\lambda\beta \intR (x - \x)^2 f^{\alpha + 1}(x, t)\, \d x + \sigma^2\rho,
\end{aligned}
\)
which presents strong similarities with the energetic bound initially studied in~\cite{toscani2024condensation}.
\end{remark}

In the same fashion as for the estimates~\eqref{eq:ptrick} and~\eqref{eq:Wbound2}---and using the same notation for the constants $c_\alpha$ and $d_\alpha$---we can control the integral comprising $f^{\alpha + 1}$ as
\(
\intR (x - \x)^2 H^{2\alpha + 2}(x) f^{\alpha + 1}(x,t)\, \d x
	\ge \frac{\biggl(\intR H^2(x)f(x,t)\, \d x\biggr)^{\frac{3\alpha}{2}}}{(c_\alpha + d_\alpha^{-2})^{\frac{3\alpha}{2}}U^{\frac{\alpha - 2}{2}}(t)},
\)
for all $\alpha > 2$. Therefore, the estimates~\eqref{eq:Hevolution} and~\eqref{eq:H2evolution} lead to
\(\label{eq:Uevolutionfinal}
	\dt U(t)
	\le - 2\lambda(\eta + 1) U(t) - \frac{\Theta_{\alpha,\eta}}{U^{\frac{\alpha - 2}{2}}(t)}+  \Xi_{\alpha,\eta},
\)
where we define
\[
\Theta_{\alpha,\eta} \coloneqq
	\biggl(\frac{2\lambda(\eta + 1)\beta}{c_\alpha + d_\alpha^{-2}} \intR H^2(x)f(x,0)\, \d x\biggr)^{\frac{3\alpha}{2}},
	\quad
\Xi_{\alpha,\eta} \coloneqq
\sigma^2(\eta + 1)(2\eta + 1)\intR H^4(x) f(x,0)\, \d x.
\]
Arguing like for the estimate~\eqref{eq:Vaftertrick}, we can now disregard the nonpositive term $- 2\lambda(\eta + 1) U(t)$, so that we are left with considering
\(\label{eq:Uevolutionsimplified}
\dt U(t) \le - \frac{\Theta_{\alpha,\eta}}{U^{\frac{\alpha - 2}{2}}(t)}+  \Xi_{\alpha,\eta}.
\)
We can observe that the right-hand side of~\eqref{eq:Uevolutionsimplified} at the time~$t=0$ is negative, provided that~$U(0)$ is small enough. On the other hand, if we introduce the  functional
\(\label{eq:Tdef}
\mathcal T(t) \coloneqq \frac{1}{\rho}\intR (x - \x)^{2\eta + 2} f(x,t)\, \d x = \frac{1}{\rho} U(t),
\)
from~\eqref{eq:Uevolutionsimplified} we get that its evolution is given by
\(\label{eq:Tevolution}
\dt \mathcal T (t)
	\le \frac{1}{\rho}\biggl(\Xi_{\alpha,\eta} - \frac{\Theta_{\alpha,\eta}}{\mathcal T^{\frac{\alpha - 2}{2}}(t)\rho^{\frac{\alpha - 2}{2}}}\biggr).
\)
From~\eqref{eq:Tevolution} it follows that if the initial mass $\rho$ satisfies
\(\label{eq:massconstraint}
\rho > \biggl(\frac{\mathcal T^{\frac{\alpha - 2}{2}}(0)\Xi_{\alpha, \eta}}{\Theta_{\alpha, \eta}}\biggr)^{\frac{2}{\alpha}},
\)
then the right-hand side of~\eqref{eq:Tevolution} is negative at time~$t=0$. Therefore, if either the initial value of the energy-type functional~$U(0)$ is sufficiently small or the initial mass is sufficiently large, we have that
\[
\dt U(t) \le -\frac{\Theta_{\alpha, \eta} - \Xi_{\alpha,\eta}U^{\frac{\alpha-2}{2}}(0)}{U^{\frac{\alpha-2}{2}}(t)},
\]
which implies
\[
U^{\frac{\alpha}{2}}(t) \le U^{\frac{\alpha}{2}}(0) - \frac{\alpha}{2}
\Bigl[\Theta_{\alpha, \eta} - \Xi_{\alpha,\eta}U^{\frac{\alpha-2}{2}}(0)\Bigr] t,
\]
so that at time
\[
t^* = \frac{2U^{\frac{\alpha}{2}}(0)}{\alpha\Bigl[\Theta_{\alpha, \eta} - \Xi_{\alpha,\eta}U^{\frac{\alpha-2}{2}}(0)\Bigr]} 
\]
we have $U(t^*) = 0$ for all $\alpha > 2$ and $\eta < 1/4$. Notice that $t^* > 0$ provided that
\[
U(0)< \biggl[\frac{\Theta_{\alpha, \eta}}{\Xi_{\alpha,\eta}}\biggr]^{\frac{\alpha-2}{2}}.
\]
These computations allow us to conclude that the solution to~\eqref{eq:sol2} loses regularity in the weighted space $L^2_w$ in finite time.

\section{A marginal-based superlinear CBO method}\label{sec:marginal}

Extending the one-dimensional results presented in the previous section to higher dimensions is challenging, as is the practical use of the system of SDEs~\eqref{eq:general-sde}. 
This difficulty is primarily due to the prohibitive computational cost associated with reconstructing the probability density via multidimensional histograms, which suffers from the curse of dimensionality. For this reason, we propose an alternative formulation that preserves the key features of the one-dimensional system while significantly reducing the computational burden in higher-dimensional simulations.

Motivated by the superior performances of anisotropic diffusion highlighted in~\cite{carrillo2021consensus}, that refined the initial general approach~\cite{pinnau2017consensus}, and by the need to derive an algorithm that scales linearly with the number of particles, we consider a multidimensional version of system~\eqref{eq:general-sde} in which, instead of taking the whole empirical measure $f_N$ in the drift term, we rather consider a mollification of the \emph{directional marginal}
\[
f_{\epsilon,N}^j(x,t) = \frac1N \underbrace{\intR\dots\intR}_{d-1 \text{ times}}\ \sum_{j = 1}^N \psi_\epsilon(x - x_1(t))\dots \psi_\epsilon(x - x_d(t))\, \d x_1 \dots \d x_{j - 1}  \d x_{j + 1} \dots \d x_d,
\]
where $\epsilon>0$. 

This is implemented by considering a novel system of SDEs written component-wise in the form
\(\label{eq:multid-SDE}
\d X_i^j = -\lambda (X_i^j - (\xxa)^j)(1 + \beta K(X_i^j)f_{\epsilon,N}^{j}(X_i^j,t)^\alpha)\d t + \sigma  (H(X_i^j))\d B_i^j,
\)
where $i = 1,\ldots, N$ and $j = 1, \ldots, d$. An Euler-Maruyama scheme applied to system~\eqref{eq:multid-SDE} essentially implements a one-dimensional algorithm of type~\eqref{eq:SLCBO1d} where coupling occurs through the marginals, so that we retain the advantages given by the one-dimensional implementation avoiding the full multi-dimensional reconstruction of the histogram density. Note that, remarkably, for $\beta=0$ the system still yields the classical CBO dynamic for $H(X_i^j)=X_i^j - (\xxa)^j$. Furthermore, except for the coupling given by the consensus point $\xxa$, the SDE for each marginal is decoupled from the other.

\subsection{Mean-field limit}
We now extend the formal mean-field derivation of Section~\ref{sec:model} to the marginal-based consensus dynamics introduced in equation~\eqref{eq:multid-SDE}. 

We now extend the mean-field derivation to the marginal-based formulation, under the same assumptions \ref{hypothesis:1}--\ref{hypothesis:4}, applied componentwise to the directional marginals $f^j(x^j,t)$ and the component of the consensus point $X_\gamma(t) \in \mathbb{R}^d$ defined as 
\[
X_\gamma^j(t) = \frac{\int_{\mathbb{R}^d} x^j e^{-\gamma F(x)} f(x,t) dx}{\int_{\mathbb{R}^d} e^{-\gamma F(x)} f(x,t) dx} \quad \text{with } f(x,t) = \prod_{j=1}^d f^j(x^j,t).
\]
Let $\varphi_j \in C_c^\infty(\mathbb{R})$ be a smooth test function. By applying Itô's formula to $\varphi_j(X_i^j(t))$, we obtain:
\begin{align*}
d\varphi_j(X_i^j(t)) &= \varphi_j'(X_i^j(t))\, dX_i^j(t) + \frac{1}{2} \varphi_j''(X_i^j(t))\, (dX_i^j)^2 \\
&= -\lambda \varphi_j'(X_i^j)\left(X_i^j - X_\gamma^j\right)\left(1 + \beta K(X_i^j)(f_{\varepsilon,N}^j(X_i^j))^\alpha \right)dt + \frac{\sigma^2}{2} H^2(X_i^j) \varphi_j''(X_i^j) dt \\
&\quad + \sigma H(X_i^j)\varphi_j'(X_i^j)\, dB_i^j.
\end{align*}
Averaging over the particles and using the empirical measure $f_{\varepsilon,N}^j(x,t)$, we obtain:
\begin{align*}
\frac{d}{dt} \left\langle f_{\varepsilon,N}^j, \varphi_j \right\rangle
&= \int_{\mathbb{R}} \left[ -\lambda (x - X_\gamma^j)\left(1 + \beta K(x)(f_{\varepsilon,N}^j(x))^\alpha\right) \varphi_j'(x)
+ \frac{\sigma^2}{2} H^2(x)\, \varphi_j''(x) \right] f_{\varepsilon,N}^j(x)\, dx.
\end{align*}
Passing formally to the limit as $N \to \infty$ and $\epsilon \to 0$, under assumptions \ref{hypothesis:1}--\ref{hypothesis:4}, we obtain the weak form of the limiting PDE:
\[
\frac{d}{dt} \langle f^j(t), \varphi_j \rangle =
\int_{\mathbb{R}} \left[ -\lambda (x - X_\gamma^j)\left(1 + \beta K(x)(f^j(x,t))^\alpha\right)\varphi_j'(x)
+ \frac{\sigma^2}{2} H^2(x)\, \varphi_j''(x) \right] f^j(x,t)\, dx.
\]
Hence, the corresponding strong form reads:
\begin{equation}\label{eq:meanfield_marginal}
\partial_t f^j(x,t) = \lambda \partial_{x} \left( (x - X_\gamma^j)\left(1 + \beta K(x)(f^j(x,t))^\alpha\right) f^j(x,t) \right) + \frac{\sigma^2}{2} \partial_{x x} \left( H^2(x) f^j(x,t) \right).
\end{equation}
Although equation~\eqref{eq:meanfield_marginal} is one-dimensional for each component $j$, the dynamics are globally coupled through the common consensus point $X_\gamma(t)$.

\subsection{Finite time blow up of the marginals}

Similarly to the one-dimensional case, let us consider the second moment functional for each marginal:
\begin{equation}
    \mathcal{V}^j(t) \coloneqq \int_{\mathbb{R}} (x - x_*^j)^2 f^j(x,t) dx, 
\end{equation}
where $x_*^j$ is the $j$-component of the unique global minimum of $\mathcal{F}(x)$.

The evolution of $\mathcal{V}^j(t)$ satisfies :
\begin{align*}
    \frac{d}{dt} \mathcal{V}^j(t) &= -2 \mathcal{V}^j(t) + 2 (X_\gamma^j - x_*^j) \int_{\mathbb{R}} (x - x_*^j) f^j(x,t) dx \\
    &\quad + 2 \sigma \int_{\mathbb{R}} H^2(x) f^j(x,t) dx - 2 \beta \int_{\mathbb{R}} (x - x_*^j)^2 K(x) (f^j)^{\alpha+1}(x,t) dx \\
    &\quad + 2 \beta (X_\gamma^j - x_*^j) \int_{\mathbb{R}} (x - x_*^j) K(x) (f^j)^{\alpha+1}(x,t) dx.
\end{align*}

Under the Laplace principle, $X_\gamma^j \to x_*^j$ as $\gamma \to \infty$, and the last two terms can be handled as perturbations. The key nonlinear damping term remains:
\begin{equation}
    A^j(t) \coloneqq \int_{\mathbb{R}} (x - x_*^j)^2 K(x) (f^j)^{\alpha+1}(x,t) dx,
\end{equation}
which can be bounded below in terms of $\mathcal{V}^j(t)$ and $\|f^j\|_{L^2}$ using the interpolation inequality:
\begin{equation} \label{eq:interpolation}
    A^j(t) \geq C_{\text{int}} \frac{(\mathcal{V}^j(t))^{(\alpha+1)/2}}{\|f^j\|_{L^2}^\alpha}.
\end{equation}
The estimate \eqref{eq:interpolation} follows from a weighted interpolation inequality between the $L^2$ norm and the second-order moment of $f$ (see~\cite{toscani2024condensation}). Following the same arguments of Section \ref{sec:meanfield}, this yields:
\begin{equation}
    \frac{d}{dt} \mathcal{V}^j(t) \leq -c_1 \mathcal{V}^j(t) - c_2 (\mathcal{V}^j(t))^{\frac{\alpha - 2}{2}},
\end{equation}
under appropriate bounds on $\|f^j\|_{L^2}$ and smallness of $C_\gamma$. This implies finite-time loss of the $L^2$ regularity  of the marginal  when $\alpha > 2$ independently on the dimension.  % oncentration (vanishing of $\mathcal{V}^j$) when $\alpha > 2$.

\subsubsection*{Critical mass and weighted energy decay}

Assuming $K(x) = H^{2\alpha}(x)$ and $H(x) = |x - x_*^j|^\eta$ with $\eta < 1/2$, we define the weighted energy functional:
\begin{equation}
    U^j(t) \coloneqq \int_{\mathbb{R}} (x - x_*^j)^2 H^2(x) f^j(x,t) dx = \int_{\mathbb{R}} |x - x_*^j|^{2 + 2\eta} f^j(x,t) dx.
\end{equation}

Then $U^j(t)$ satisfies the differential inequality (cf. equation~\eqref{eq:Uevolutionsimplified}):
\begin{equation}
    \frac{d}{dt} U^j(t) \leq - \Theta_{\alpha,\eta} (U^j(t))^{\frac{\alpha - 2}{2}} + \Xi_{\alpha,\eta},
\end{equation}
for constants $\Theta_{\alpha,\eta}, \Xi_{\alpha,\eta} > 0$ depending on the initial mass and $\sigma$. Hence, if $U^j(0)$ is small enough or the mass is sufficiently large, one observes loss of weighted $L^2$ regularity in finite time.

Although the marginal-based formulation significantly simplifies the multidimensional setting, it still retains the essential features of the original dynamics, including the emergence of finite-time singularities.

\section{Numerical tests}\label{sec:numerics}
In this section, we present a series of numerical experiments that 
support our analytical findings. These are divided into two parts: the first 
one focuses on one-dimensional simulations, while the second one
introduces a modified scheme tailored to higher dimensions. This 
multidimensional variant, based on marginalized interactions, 
demonstrates the potential to outperform the standard CBO method when 
implemented within the same computational framework.

\subsection{One-dimensional examples}\label{sec:1d}
In the one-dimensional tests, we consider the following system of SDEs
\(\label{eq:generic-SDE}
\d X_i = -\lambda (X_i - \xxa)(1 + \beta K(X_i)f^\alpha_{\epsilon,N}(X_i,t))\d 
t + \sigma  H(X_i)dB_i,
\)
where $X_i \in [-L, L]^d$, and $i = 1,\ldots, N$ and $\epsilon>0$ is a numerical parameter characterizing the density reconstruction which depends on the mesh size and will be made explicit later. Equation~\eqref{eq:generic-SDE} is
discretized with the Euler-Maruyama forward scheme such that, the $i$-th trajectory at the iterate $n+1$ is given by
\(
X_i^{n+1} = X_i^{n} + \lambda \Delta t_{n} (\xxa^{n} - X_i^{n})(1 + \beta K(X_i^{n}) f_{\varepsilon,N}^\alpha(X_i^{n}))
	+ \sigma\sqrt{\Delta t_{n}}  H(X_i^{n})\xi_i^n,
	\label{eq:SLCBO1d}
\)
where $\xi_i^n$ is distributed as $N(0,1)$.

The time step $\Delta t_n$ is chosen adaptively to ensure that the numerical trajectories remain confined within the search space. In particular, due to the presence of the local density $f_{\epsilon,N}^\alpha$ in the drift term, which may become large in regions of particle concentration, the dynamics may become strongly amplified. To prevent overshooting and preserve numerical stability, we impose the condition
\(\label{eq:CFL}
\Delta t_n \le \min_i \frac{1}{\beta 
K(X_i^n)f_{\epsilon,N}^\alpha(X_i^n)},
\)
which balances the intensity of the interaction with a corresponding reduction in the time increment. This choice becomes crucial in condensation regimes, where particles cluster tightly and the drift term grows rapidly.

For the reconstruction of $f_{\epsilon}^N(X_i^n)$, we employ a histogram-based approach using the mollifier $\psi_\epsilon(y) = 1/\epsilon \chi(|y| \leq \epsilon/2)$ with $\epsilon = \Delta x$. As observed in~\cite{calzola2024emergence}, the number of bins used in the histogram plays a critical role: a sufficiently fine discretization is required to correctly capture the formation of the steady state, especially in the subcritical regime. In all simulations, we fix $\lambda = 1$, $\sigma^2=1/4$ and adopt a staggered mesh to ensure that the point of symmetry $\x$ does not coincide with any grid node, avoiding numerical artifacts due to the solution's singularity at that point.

\subsubsection[\texorpdfstring{Case $H(x) = (x - \x)^\eta$ and $K(x) = H^\alpha(x)$}{Case H(x) = (x - x̄)²ⁿ and K(x) = Hᵃ(x)}]{Case $\bm{H(x) = (x - \x)^{\eta}}$ and $\bm{K(x) = H^\alpha(x)}$}

Each test aims to highlight the behavior of the solution to the PDE
\(\label{eq:pdeHK}
\pd f(x,t) = \pd[x]\Bigl(\lambda\bigl[(x - \xa(t))f(x,t)(1 + \beta (x - \xa(t))^{2\eta\alpha} f^\alpha(x,t)\bigr] + \frac{\sigma^2}{2}\pd[x] \bigl[(x - \xa(t))^{2\eta}f(x,t)\bigr]\Bigr),
\)
either in the case $\xa(t) \equiv \x$ and $\xa(t)$ defined as in~\eqref{eq:xadef}.

\subsubsection*{Consistency}
\begin{figure}[htbp]
	\centering
\begin{tblr}{}
\includegraphics[width=0.4\textwidth]{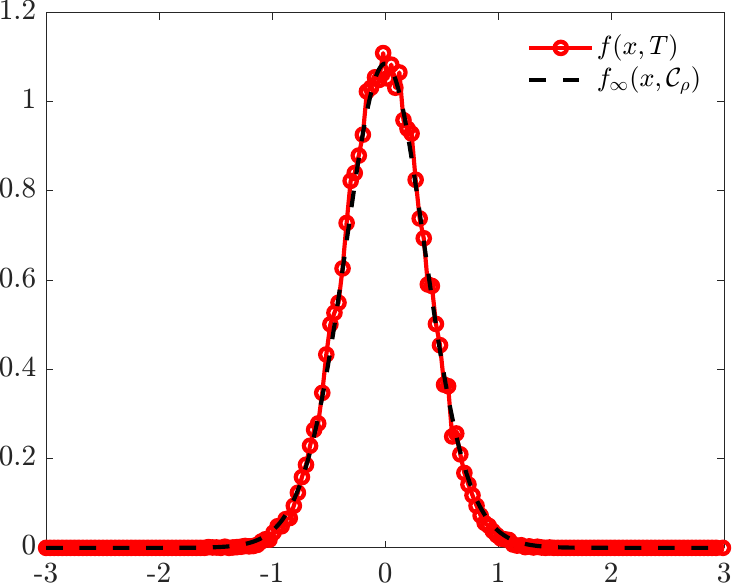} &
\includegraphics[width=0.4\textwidth]{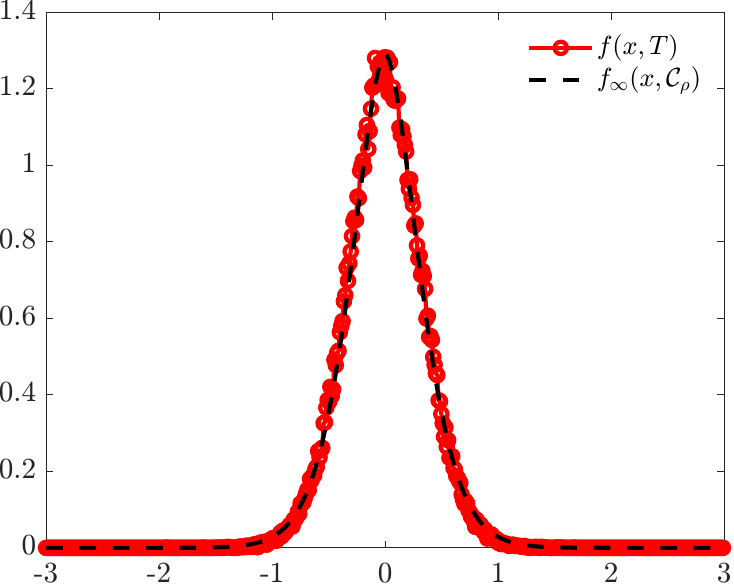} \\
\includegraphics[width=0.4\textwidth]{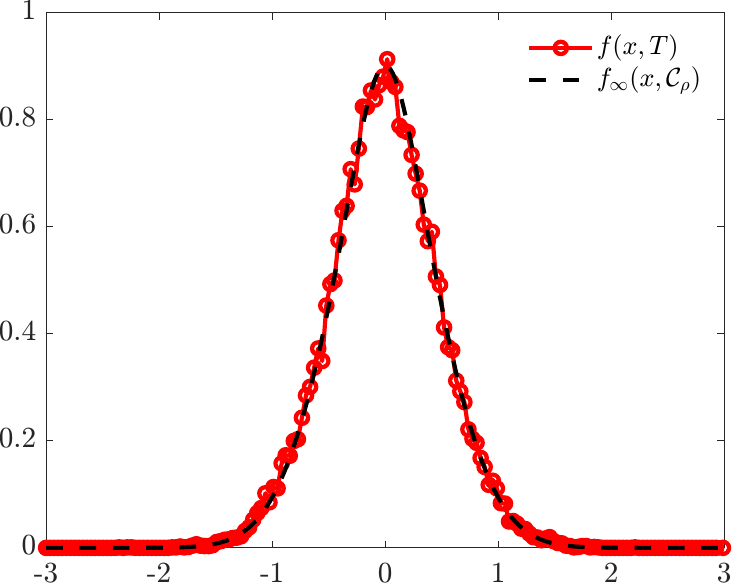} &
\includegraphics[width=0.4\textwidth]{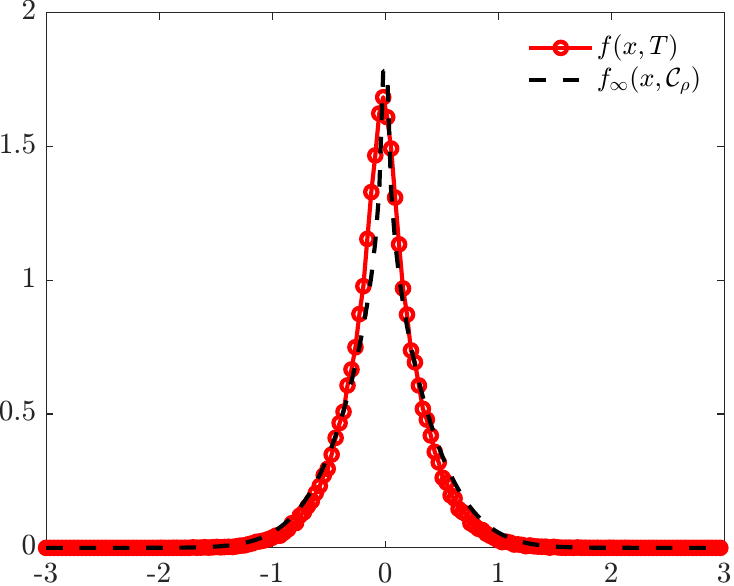} \\
\includegraphics[width=0.4\textwidth]{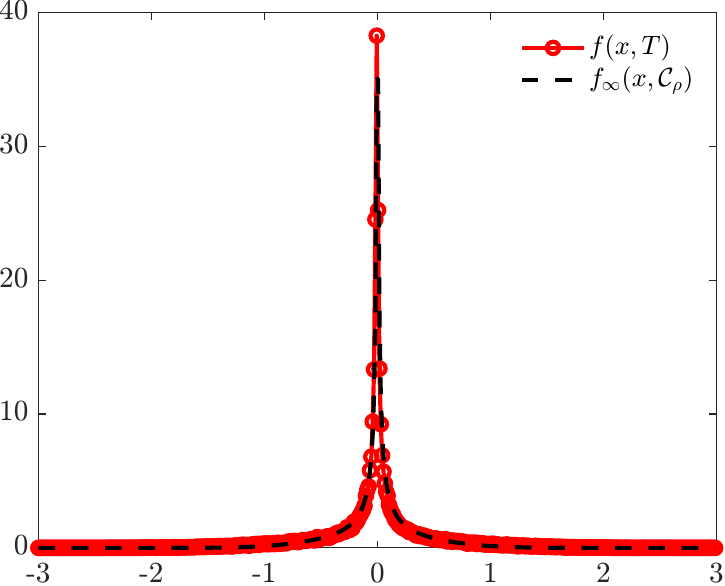} &
\includegraphics[width=0.4\textwidth]{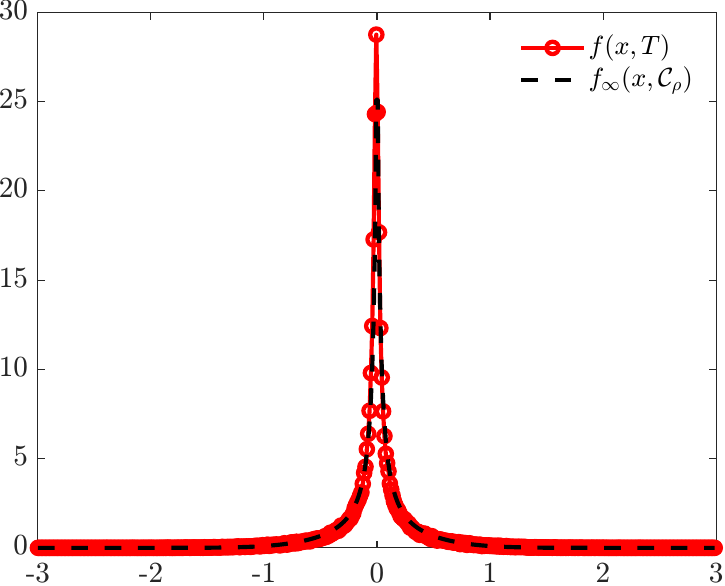} 
\end{tblr}
\caption{Steady states. Left, $\eta = 0$, right, $\eta = 1/8$. Top: $\alpha = 0.25$ and subcritical mass, middle: $\alpha = 2.25$ and subcritical mass, bottom: $\alpha = 2.25$ and supercritical mass.}
\label{fig:consistency}
\end{figure}

In the first test, we simulate algorithm~\eqref{eq:generic-SDE} with $N = 10^5$ particles, fixing $\xxa(t) = \xx = 0$ for all times and choosing $H(X) = \abs{X - \xx}^{\eta}$ and $K(X) = H^\alpha(X)$. This test shows the consistency of the large-time behavior of the numerical solution to equation~\eqref{eq:pdeHK} when $N \gg 1$, both in the subcritical and supercritical case. We compute the critical mass $\bar \rho$ as the integral
\[
\bar \rho = \intR \finf(x,\bar{\mathcal C}_\rho)\, \d x, \qquad \bar{\mathcal C}_\rho = \frac{1}{\beta^\alpha}.
\]
We test for the cases $\alpha = 0.25$ and $\alpha = 2.25$ and $\eta = 0$ and $\eta = 0.125$. Figure~\ref{fig:consistency} shows accordance between the analytical steady state and the numerical solution in the different regimes tested.

\subsubsection*{Role of the number of bins}

\begin{figure}[htbp]
	\centering
\begin{tblr}{}
\includegraphics[width=0.4\textwidth]{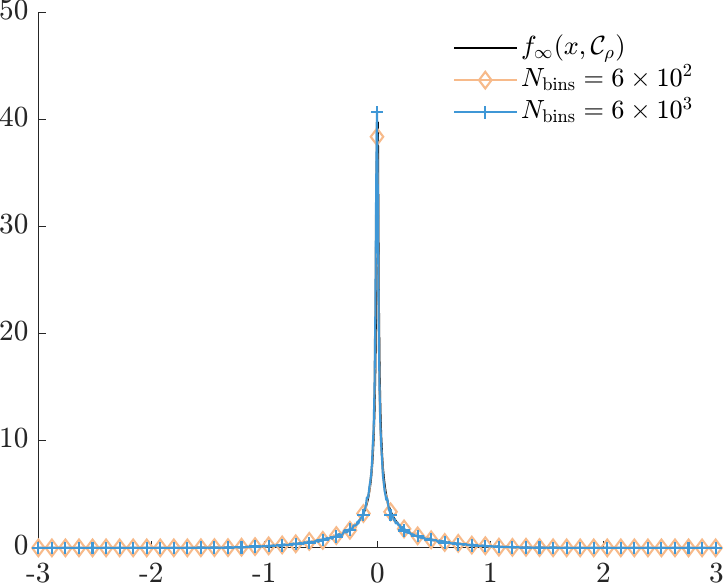} &
\includegraphics[width=0.4\textwidth]{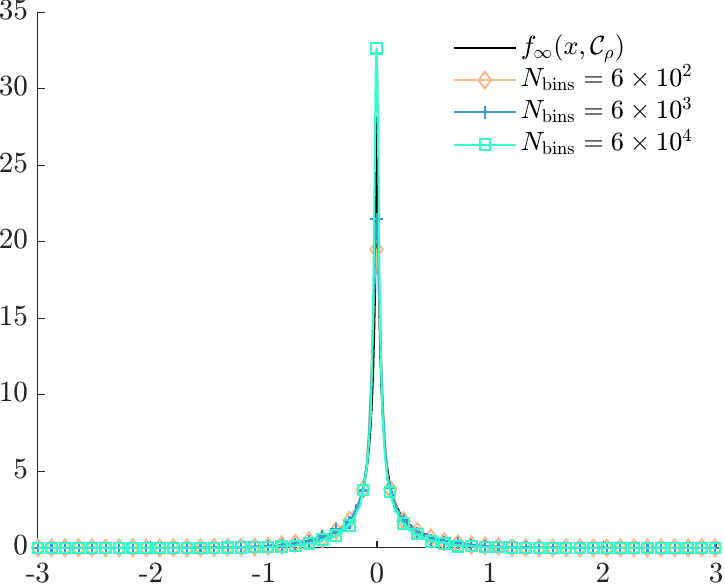} \\
\includegraphics[width=0.4\textwidth]{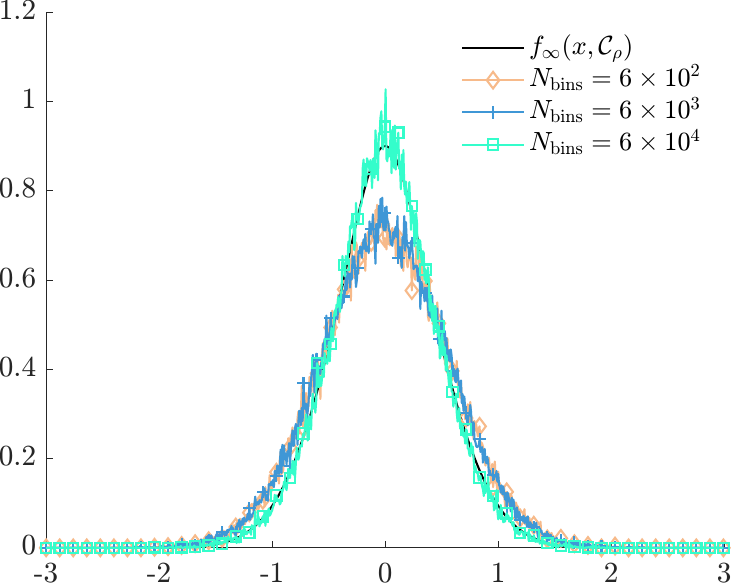} &
\includegraphics[width=0.4\textwidth]{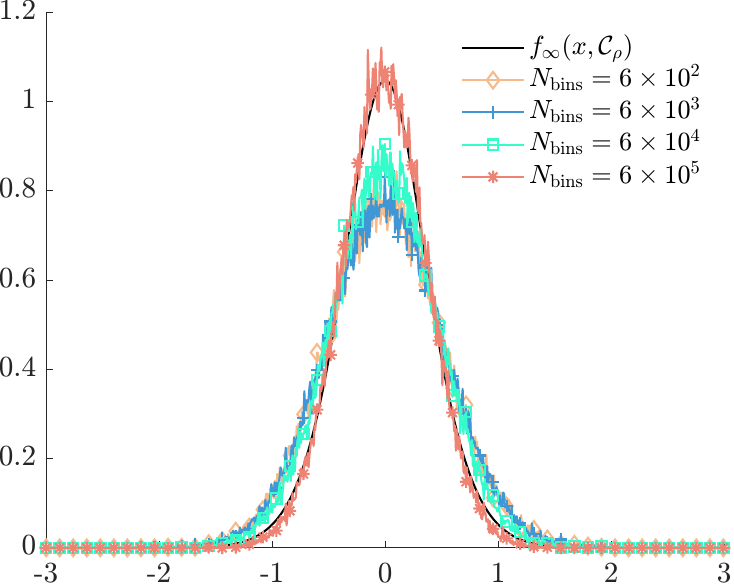} \\
\includegraphics[width=0.4\textwidth]{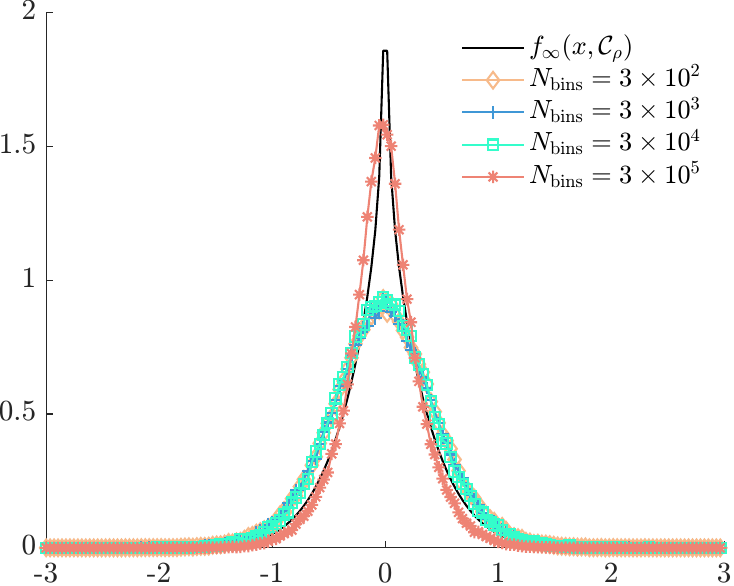} &
\includegraphics[width=0.4\textwidth]{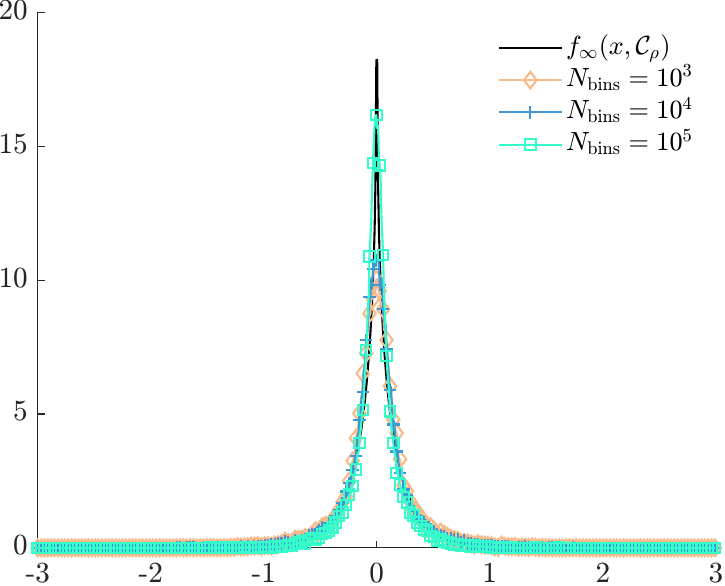} 
\end{tblr}
\caption{Steady states. Top: $\alpha = 2.25$ and supercritical mass. Left: $\eta = 0$, right: $\eta = 1/8$. Middle: $\eta = 0$ and subcritical mass. Left: $\alpha = 2.25$, right $\alpha = 0.5$. Bottom: $\alpha = 1.5$ and $\eta = 1/8$. Left: subcritical mass, right: supercritical mass.}
\label{fig:bins}
\end{figure}

We show the dependence of the number of bins on the accuracy of reconstruction of the steady state. We use again $N = 10^5$ particles to perform the simulations, testing for different values for $\alpha$, both in the case $\eta = 0$ and $\eta = 0.125$. This test highlights that the number of bins required to accurately reproduce the steady state can be high (up to $\mathcal O(10^5)$ bins in the simulations shown in Figure~\ref{fig:bins}), and that the case $\eta > 0$, which implies the presence of a blow-up at $x = \x$ regardless of the initial mass, presents more challenges with respect to the case $\eta = 0$. This is consistent with the findings firstly reported in~\cite{calzola2024emergence}.

\subsubsection*{Performance in optimization}\label{sec:511}

We consider the objective function $\F(x) = x^2$ on the domain $[-3, 3]$ and we take $N = 10^5$ particles to find its minimizer $\x = 0$. We set $\alpha = 2.25$ and consider a supercritical mass for both $\eta = 0$ and $\eta = 0.125$. In particular, we fix $\rho = 12$ in both cases. We compute $\xxa(t)$ to be
\[
\xxa(t) = \frac{\int_{-3}^{3} xe^{-\gamma\F(x)} f(x,t) \, \d x}
{\int_{-3}^{3} e^{-\gamma\F(x)} f(x,t)\, \d x}, \qquad \gamma = 50,
\]
using $N_{\mathrm{bins}} = 201$ bins for the reconstruction of $f(x,t)$ at each iteration.

\begin{figure}
	\centering
\begin{tblr}{}
	\includegraphics[width=0.4\textwidth, trim=2em 1em 3em 1em]{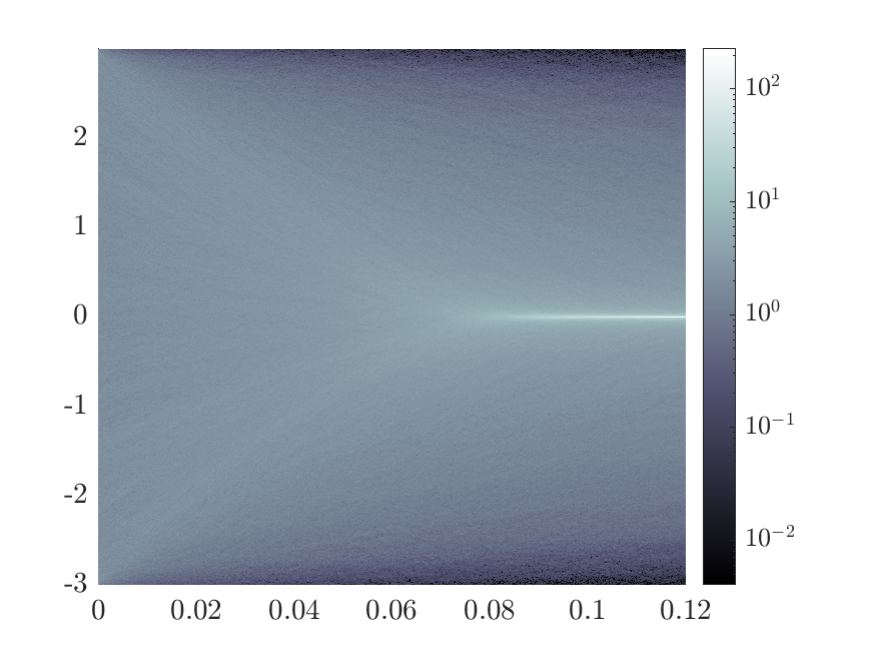} &
	\includegraphics[width=0.4\textwidth, trim=2em 1em 3em 1em]{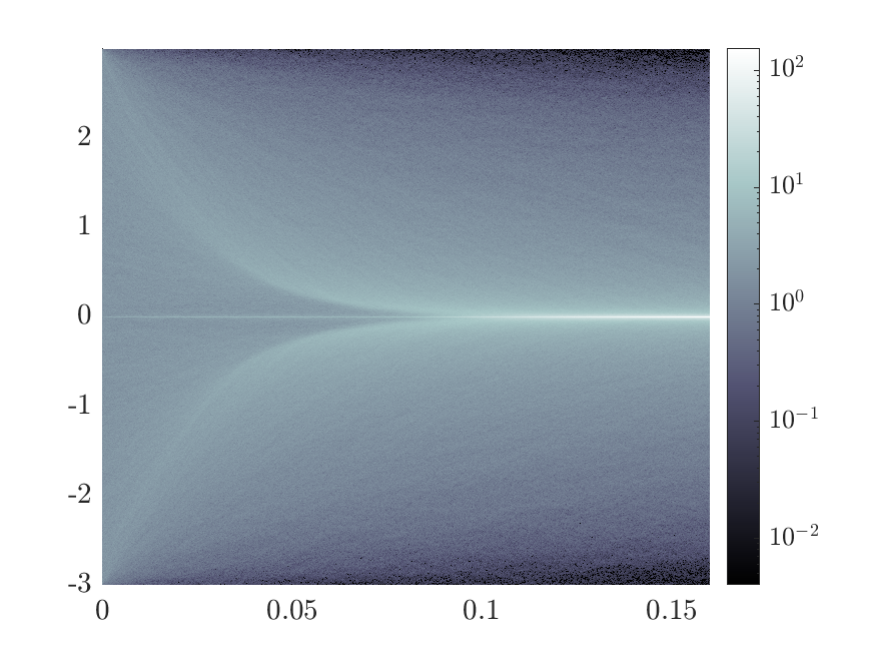}
\end{tblr}
\caption{Evolution of the solution $f(x,t)$ to equation~\eqref{eq:pdeHK} applied to the minimization of the parabola $\F(x) = x^2$ on $[-3, 3]$. The initial mass is $\rho = 12$ (a supercritical value). Left: the case $\eta = 0$; right: the case $\eta = 0.125$. In both cases we see the sharp loss of regularity with associated blow-up at the point $x = \x = 0$.}
\label{fig:f-evo-eta}
\end{figure}

We report in Figure~\ref{fig:f-evo-eta} the evolution of the numerical solution to~\eqref{eq:pdeHK} with time-dependent global best~$\xxa(t)$. As expected from the theoretical analysis, the solution quickly loses its regularity, with the formation of blow-ups before $t = 0.2$ for both $\eta = 0$ and $\eta = 0.125$. However, the mass is still relatively sparse, meaning that the energy of the solution, which in this case coincides with the functional $\V(t)$, is still relatively large. On the other hand, condition~\eqref{eq:CFL} on the time-discretization implies that the number of iterations needed in order to significantly extend the time horizon after the blow-up formation would be very large.

\subsubsection[\texorpdfstring{Case $H(x) = |x - \xa(t)|$ and $K(x) = 1$}{Case H(x) = |x - xᵧ(t)| and K(x) = 1}]{Case $\bm{H(x) = |x - \xa(t)|}$ and $\bm{K(x) = 1$}}

Motivated by our previous test, we focus on a different choice for the coefficients $K(x)$ and $H(x)$, i.e.,
\(\label{eq:pde21}
\pd f(x,t) = \pd[x]\Bigl(\bigl[\lambda(x - \xa(t))f(x,t)(1 + \beta f^\alpha(x,t)\bigr] + \frac{\sigma^2}{2}\pd[x] \bigl[(x - \xa(t))^2 f(x,t)\bigr]\Bigr),
\)
while maintaining the same definition for $\xa(t)$.

\subsubsection*{Performance in optimization}

\begin{figure}[htbp]
	\centering
	\includegraphics[width=0.5\textwidth]{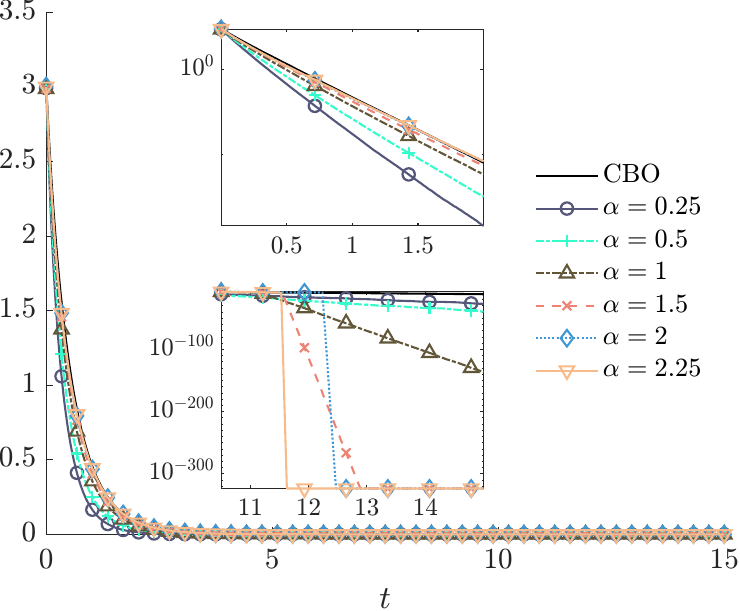}
\caption{Evolution in time of energy functional $\V(t)$ for 
equation~\eqref{eq:pde21} for different choices of~$\alpha$. We see that the 
initial decay is faster for small values of~$\alpha$, but as $\alpha \ge 1$ and 
approaching the threshold value $\bar \alpha = 2$, the energy of the solution 
decays to zero with sharper and sharper transitions. Moreover, the value at 
time $T = 15$ for $\V(T)$ for all values of~$\alpha$ is several orders of 
magnitudes smaller than the reference CBO value, corresponding to $\alpha = 
\beta = 0$. }
\label{fig:V-decays}
\end{figure}

We simulate equation~\eqref{eq:pde21} using the scheme~\eqref{eq:generic-SDE} with $N = 10^5$ particles until the final (physical) time~$T = 15$, for various values of $\alpha$. In particular, we refine the condition on the time step to be
\[
\Delta t_n \le \min\Bigl(\min_i \frac{1}{\beta 
K(X_i^n)f_{\epsilon,N}^\alpha(X_i^n)}, 0.005\Bigr).
\]
This ensures a fair comparison with the standard CBO implementation, where a fixed time-step is usually employed.
In Table~\ref{tab:niter} we summarize the number of iterations required by the 
simulation in this setting: the larger the value of~$\alpha$, the higher the 
number of iterations, due to adaptive time discretization, as expected, 
especially when $\alpha$ exceeds the value of~$1$ and approaches the critical 
value~$\alpha = 2$, where we know that the functional $\V$ vanishes in finite 
time. The results presented in Table~\ref{tab:niter} show that, for values 
of $\alpha \le 1$, the method does not actually require an adaptive time 
discretization (the number of iterations for those cases is $15/0.005 = 3000$, 
i.e., at every time step $\Delta t_n \ge 0.005$), since the particles are not 
concentrating fast enough to cause numerical instabilities. Still, as shown in 
Figure~\ref{fig:V-decays}, the concentration is significantly faster than the 
one provided by the CBO method with the fixed time discretization $\Delta t_n = 
0.005$.

\begin{table}[htbp]
    \centering
    \begin{tblr}{
        colspec={c|cccccc},
        row{1}={c},
        row{2}={c},
        hline{1,Z}={1pt},
        hline{2}={0.5pt}
    }
    $\alpha$            & 0.25 & 0.5 & 1   & 1.5  & 2     & 2.25 \\
    Iterations          & 3000 & 3000 & 3000 & 5120 & 26940 & 63583
    \end{tblr}
    \caption{Number of iterations required to simulate equation~\eqref{eq:pde21} with $N = 10^5$ particles until $T= 15$, for various choices of $\alpha$.}
    \label{tab:niter}
\end{table}

We show in Figure~\ref{fig:V-decays} the different tendencies of the decay of the functional~$\V(t)$, along with a comparison to its decay with the standard CBO implementation (corresponding to our algorithm for $\beta = 0$). If we consider the behavior at the very beginning of the evolution, we see---top box---that values of $\alpha$ that are greater than~$1$ match the exponential decay of the standard CBO, and, more in general, the lower the value of $\alpha$, the faster the initial decay. This is consistent with the nonlinearity added in the drift term.

On the other hand, the larger the value of $\alpha$, the more abrupt the decay of the~$\V(t)$ functional---bottom box. In particular, the decay becomes appreciably faster than exponential as long as $\alpha \ge 1$, and when $\alpha$ surpasses the critical value~$\alpha = 2$ then the decrease becomes almost instantaneous, dropping almost~$300$ order of magnitudes in hundredths of physical-time units.

Interestingly, we observe that the subcritical regime (e.g., $\alpha = 0.25$) not only maintains strong localization properties, but also mitigates the numerical instabilities encountered in the supercritical case ($\alpha > 1$). In fact, for $\alpha > 1$, the strong concentrating effects tend to produce overly sharp peaks and finite-time blow-up, leading to poor resolution, sensitivity to the bin size, and premature convergence toward local minima. By contrast, the subcritical choice appears to stabilize the dynamics while preserving the beneficial features of consensus formation, thus offering a compelling compromise between concentration and robustness. 

\begin{figure}[htbp]
	\centering
	\begin{tblr}{}
		\includegraphics[width=0.45\textwidth, trim=2em 1em 2em 1em]{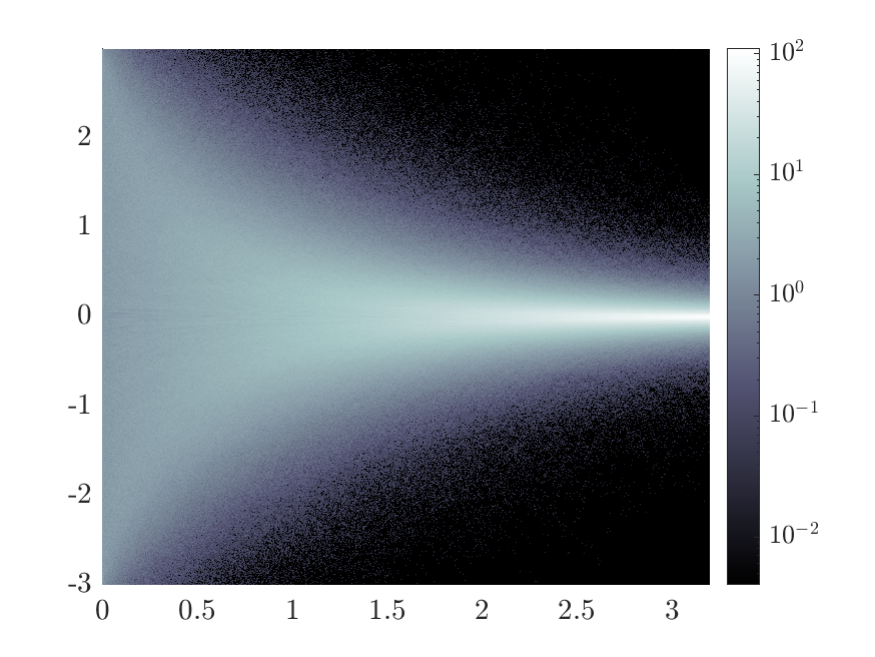} &
		\includegraphics[width=0.45\textwidth, trim=2em 1em 2em 1em]{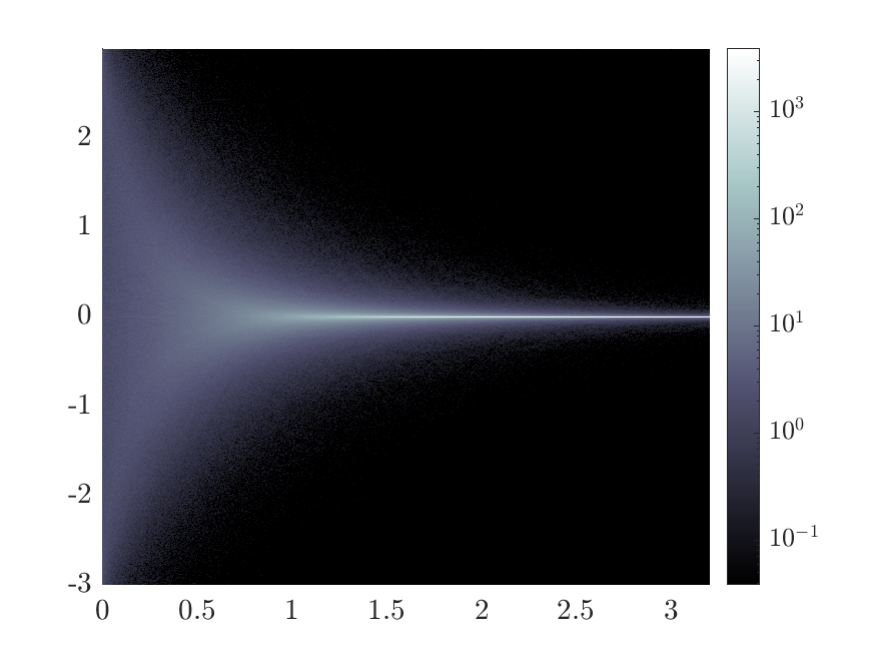}
	\end{tblr}
	\caption{Evolution of the solution $f(x,t)$ to equation~\eqref{eq:pde21} 
	applied to the minimization of the parabola $\F(x) = x^2$ on $[-3, 3]$. The 
	initial mass is $\rho = 12$. Left: standard CBO algorithm (which 
	corresponds to choosing $\beta = 0$). Right: $\alpha = 0.25$ and 
	$\beta = 1$. We see that standard CBO converges to a Dirac's Delta significantly slower opposed to its nonlinear counterpart.}
	\label{fig:f-evo-our}
\end{figure}

As a last example, we show in Figure~\ref{fig:f-evo-our} a visual comparison of the different profiles of the solution to the standard CBO model and equation~\eqref{eq:pde21} with $\alpha = 0.25$, respectively, when considering the minimization of the parabola $\F(x) = x^2$. In both cases we considered an initial mass $\rho = 12$, which is supercritical for $\alpha = 0.25$ in our setting.  We see that the additional nonlinearity is significantly more effective in concentrating the mass of the solution, all else being equal.

These results suggest that this novel approach can be effective also in higher dimensions.

\subsection{Multi-dimensional case}
As already discussed, in the multi-dimensional case our approach is based 
on the component-wise formulation introduced in Section \ref{sec:marginal}. A 
numerical scheme is then obtained as direct generalization of the 
one-dimensional case by direct application of the Euler-Maruyama method to the 
SDE system \eqref{eq:multid-SDE}. We report in Algorithm 1 the details of the 
method. 

We report in Table~\ref{tab:multi-d} the comparison of the results between 
implementation by Algorithm~\ref{alg:mscbo} of the marginal-based superlinear 
CBO~\eqref{eq:multid-SDE}, or SL-CBO in short, and classical CBO for various 
standard objective functions (see Appendix A) in dimension $d = 20$. In each 
case, with the exception of the Rosenbrock function, whose natural domain, 
which is asymmetric, has instead been mapped to $[-1.5, 3]^d$, we set the 
domain to be $[-3, 3]^d$ and we make every scenario homogeneous by applying a 
suitable scaling of variable for all functions so that their minimizer points 
actually lie in the chosen hypercube, following the standard practices of, 
e.g.,~\cite{grassi2021particle,carrillo2021consensus}. This approach is 
necessary if want to adopt the same set of parameters in all test cases. 
Furthermore, the insights gained from the one-dimensional analysis play a 
central role in selecting the value of $\alpha$. In particular, we observed 
that for $\alpha = 0.25$ the one-dimensional algorithm exhibits sharp 
concentration and improved convergence rates without compromising numerical 
stability. For this reason, we consistently choose $\alpha = 0.25$ in the 
following, along fixing $\rho = \beta = 1$ for the SL-CBO method.

\def\mystrut{\rule[-0.5ex]{0pt}{3ex}}
\AddToHook{cmd/State/after}{\mystrut}
\def\Tmax{T_{\text{max}}}
\def\nmax{n_{\text{max}}}
\begin{algorithm}[tbp]\small
\caption{Multi-dimensional CBO with superlinear drift \rule[-1.5ex]{0pt}{4.5ex}}
\label{alg:mscbo}
\begin{algorithmic}[1]
\State Set algorithmic parameters $N$, $d$, $\gamma$, $\Tmax$, $\nmax$, $\delta_{\textrm{threshold}}$
\State Set interaction parameters $\alpha$, $\beta$, $\rho$, $\lambda$, $\sigma$
\State Initialize $t = 0$, $n = 0$
\State Sample $X_i \in \R^d$, $i = 1, \ldots, N$ from the given initial distribution $f(x,0)$
\While{\rule[-0.5ex]{0pt}{3.5ex}$t < \Tmax$ \textbf{or} $n \le \nmax$}
	\State Compute $\displaystyle \xxa = \frac{1}{\sum_{i = 1}^N\, e^{-\gamma \F(X_i)}}\sum_{i = 1}^N X_i e^{-\gamma \F(X_i)}$
	\For{\rule[-0.5ex]{0pt}{3.5ex}$j = 1, \ldots, d$}
		\State Reconstruct $f_N^j(X_i^j)$ by computing the 1-D histogram of $X_i^j \in \R$
	\EndFor\mystrut
	\State Require: $\displaystyle \Delta t_n \le \min_i \frac{\strut1}{\beta K(X_i^j)[f_N^j(X_i^j)]^\alpha}$
	\State Update $X_i$ by
	\[
	\begin{split}
	X_i
		&\gets X_i - \Delta t_n\lambda\sum_{j=1}^d \bm{e}_j[(X_i - (\xxa))_j(1 + \beta K(X_i))]_j\,f_{\epsilon,N}^{j}(X_i,t)^\alpha\\
		&\hphantom{{}\gets} + \sqrt{\Delta t_n}\sigma  \sum_{j=1}^d \bm{e}_j(H(X_i))_j\xi_{i,n}^j, \quad \xi_{i,n}^j \sim \mathcal N(0, 1)
	\end{split}
	\]
	\State where $\bm{e}_j$ denotes the $j$-th standard basis vector in $\mathbb{R}^d$.
	\State Compute particle variance for the stopping criterion: $\mathrm{Var}(X_i)$
	\If{$\mathrm{Var}(X_i) \le \delta_{\textrm{threshold}}$}
		\State \textbf{break};
	\EndIf
\EndWhile\mystrut
\end{algorithmic}
\end{algorithm}

\begin{table}[htbp]
    \centering
    \includegraphics[width=\linewidth]{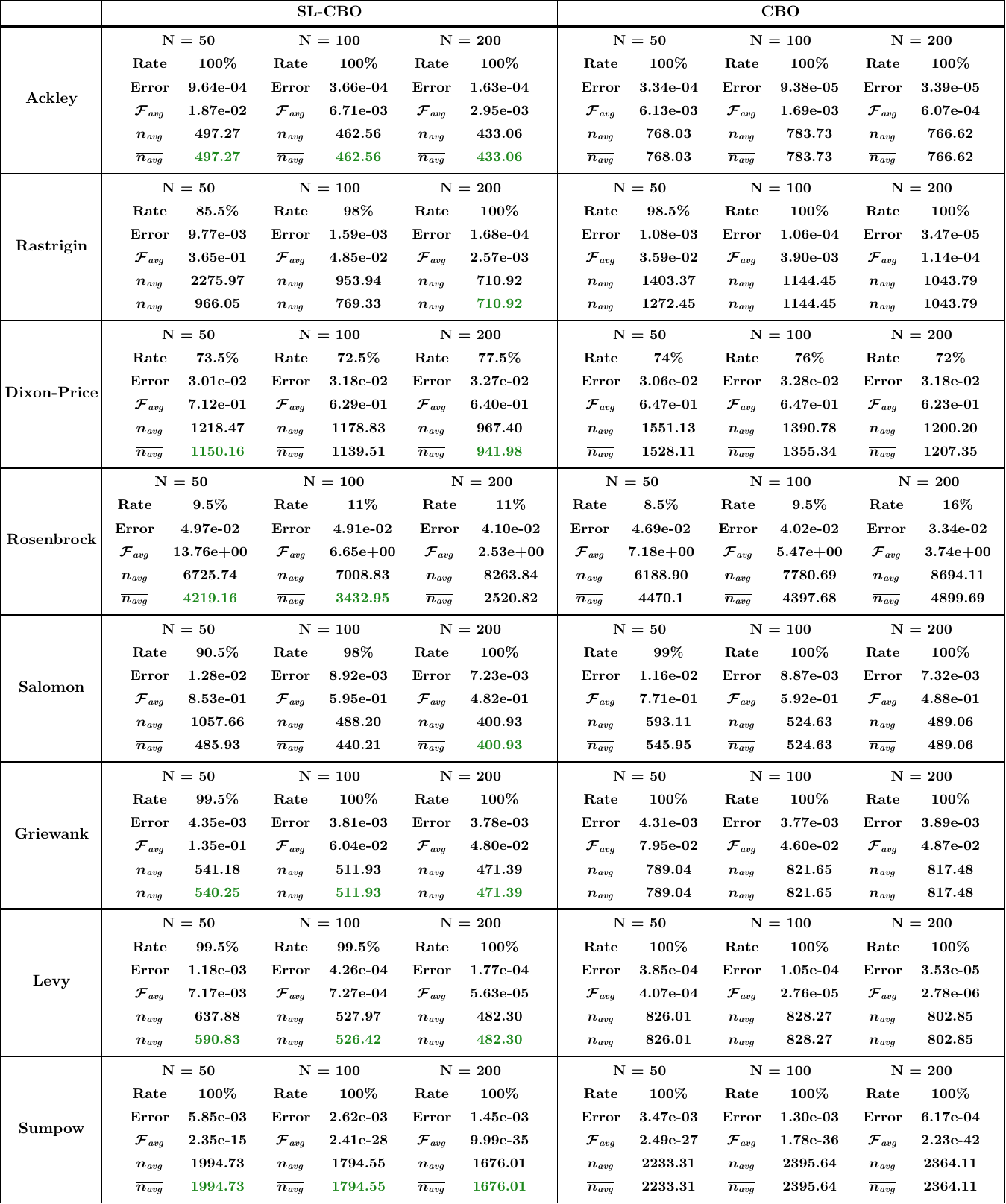}
    \caption{Comparison between SL-CBO and CBO algorithm for different test functions for global optimization (see Appendix A) and varying number of particles. All tests were performed with $d=20$, $\sigma = 5$ and $\Delta t =0.05$. We highlighted in green the instances in which the SL-CBO method returns the correct global minimizer accurately (within $0.5\%$ difference in success rate with respect to CBO) while employing less iterations than its classical counterpart.}
	\label{tab:multi-d}
\end{table}

For each objective function, we report the success rate over $N_{\text{sim}} = 200$ simulations, which is defined as the percentage of simulations such that, at the end of the run, the condition
\[
\max_j \abs{(\xa)_j - (\x)_j} < \delta_{\text{threshold}}
\]
is satisfied. We chose $\delta_{\text{threshold}} = 0.25$ as 
in~\cite{carrillo2021consensus}. Each run ends either if the particle variance 
drops below $10^{-2}$, or if the total number of iterations reaches 
$10^4$. We perform every simulation using $N = 50$, $100$ and~$200$ particles 
and fixing $\lambda = 1$ and $\sigma = 5$, again 
following~\cite{carrillo2021consensus}. This higher value for the diffusion 
parameter (in spite of the theoretical value for ensuring convergence being 
$\sigma < \sqrt2$), is chosen in order to allow the exploratory dynamics 
inherent of these methods to compensate for the low number of particles used, 
which clearly is far from satisfying $N \gg 1$. This choice is effective in 
ensuring very good numerical performance, especially for the objective 
functions that present many local minima and are oscillatory in nature, 
compared to those that are rather valley-shaped, where both our method and 
classical CBO tend to trail behind.

We report the final average $\ell^2$ error between the final value for 
$\xa$ and the actual minimizer $\x$, the average value of $\F(\xa)$ at final 
time, and both the average number of iterations in general $n_{avg}$ and the 
average number of iterations when the simulation ended successfully 
($\overline{n_{avg}}$).

The stopping condition on the particle variance and on the number of 
iterations~$n_{\max}$, the success threshold~$\delta_{\text{threshold}}$, the 
number of particles~$N$, the consensus parameter~$\gamma$, the drift 
coefficient~$\lambda$, the diffusion coefficient~$\sigma$ and the time 
delta~$\Delta t_n$ are shared between the CBO method and the SL-CBO method, 
while $\alpha$, $\beta$ and $\rho$ are set just for the SL-CBO method. The 
number of bins used for the marginal densities' reconstruction was set as 
$N_{\text{bins}} = 2001$. The role of the number of bins for optimization 
purposes is much less prominent than the one played for capturing supercritical 
steady states (cfr.\ the optimization test in Section~\ref{sec:511}, where 
$N_{\text{bins}}$ was set equal to $201$), as the performance of the SL-CBO 
method is not significantly affected by varying the number of bins, in our 
tests. Finally, for both methods, the time discretization is identical, setting 
$\Delta t_n = \delta t = 0.05$, since, as noted in the one-dimensional case, 
the subcritical regime $\alpha = 0.25$ does not require special treatment to 
ensure the numerical stability of the time discretization.

In general, it is clear that an increase in the number of particle used implies a better performance overall, improving both errors and speed of convergence.
Both our implementation and the classical CBO perform worse on functions characterized by narrow valleys or highly oscillatory landscapes, which are known to pose significant challenges for gradient-free optimization methods. We highlighted in green the cases in which our implementation outperforms the standard CBO one, with equal setting, i.e., when the number of iterations is lesser while maintaining a success rate at least within a $0.5\%$ difference.

The general trend that we may observe is that the SL-CBO is capable of performances very similar to the ones provided by the standard CBO implementation, in terms of accuracy, especially when using $N = 100$ or $N = 200$ particles. Moreover, Table~\ref{tab:multi-d} shows several cases in which the superlinear CBO method can reach the global minimizer saving a significant number of iterations ($-30\%$ on average, when successful), proving itself useful also in high-dimensional applications. 

To ease the comparison between classical CBO and SL-CBO we measured the total runtime of both a standard CBO implementation and our SL-CBO one with the built-in function \texttt{timeit} of MATLAB. We took into account Ackley and Rastrigin costs. We highlight that this test shows the ratio of the median runtime of the SL-CBO method with respect to the CBO method (so that an output greater than 1 means that the standard CBO is faster, whereas an output lesser than 1 means that the SL-CBO method is faster). 
\begin{figure}
\centering
	\includegraphics[width=0.5\linewidth]{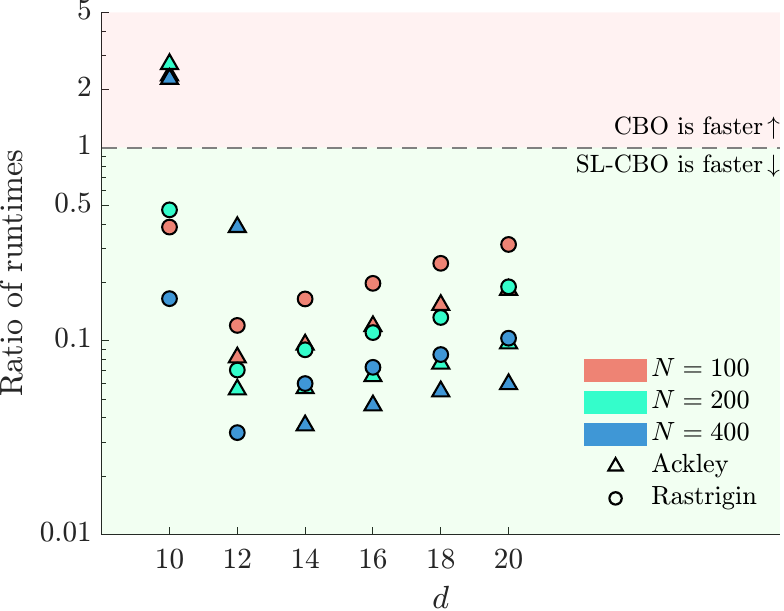}
	\caption{Comparison between SL-CBO and CBO algorithm in terms of runtime. }
	\label{fig:runtime}
\end{figure}
Each runtime was measured across 5 \texttt{timeit} calls for each combination of dimension value, number of particles and cost function. In all cases, the parameters were are the ones reported in Table \ref{tab:multi-d}. As it can be seen in Figure \ref{fig:runtime}, the reconstruction cost may outweigh its benefits, especially in lower dimensions. However, approaching $d=20$, the savings in both computing time and number of iterations suggest an advantage of SL-CBO over classical CBO dynamics.

\section{Conclusions and perspectives}
We proposed a new variant of consensus-based optimization motivated by recent extensions of the Kaniadakis--Quarati model for indistinguishable bosons. The resulting mean-field formulation leads to a nonlinear Fokker--Planck equation with superlinear drift and nonconstant diffusion, whose analysis reveals finite-time blow-up and condensation phenomena, analogously to the underlying quantum models. A major challenge is related to the multidimensionality of the resulting system which create burdens both at the analytical as well as computational problem. In one space dimension we derived the mean-field limit of the system and demonstrated the emergence of condensation effects, including the finite-time blow-up and loss of $L^2$-regularity in the solution. However convergence to global minimum remains an open problem.

To achieve a computationally efficient model we modified the original system taking inspiration by the one-dimensional analysis and adding directional superlinear effects using the marginals of the distribution function. This permits to scale the overall computational complexity linearly with the numbers of particles as in classical CBO systems.
From a computational viewpoint, preliminary numerical experiments show that the proposed scheme may offer enhanced concentration and efficiency in reaching global minimizers, especially in the presence of complex or multimodal landscapes. These results suggest that introducing nonlinear drift terms and nonuniform diffusion can potentially improve the convergence properties of CBO-type methods.

Several directions remain open for future research. First, a rigorous analysis of the convergence properties of the stochastic particle system to the mean-field limit and the global minimum is still lacking. Second, adaptive strategies for tuning the drift and diffusion terms could enhance robustness across a broader class of optimization landscapes. Finally, extending the proposed approach to constrained and multiobjective settings, or integrating it with machine learning paradigms, could further expand its applicability.

\section*{Acknowledgements}
The work of LP was partially supported by the Royal Society under the Wolfson Fellowship ``Uncertainty quantification, data-driven simulations and learning of multiscale complex systems governed by PDEs\rq\rq. LP also acknowledges the partial support by the European Union through the Future Artificial Intelligence Research (FAIR) Foundation, ``MATH4AI\rq\rq\ Project and by ICSC -- Centro Nazionale di Ricerca in High Performance Computing, Big Data and Quantum Computing, funded by European Union -- NextGenerationEU and by the Italian Ministry of University and Research (MUR) through the PRIN 2022 project (No.\ 2022KKJP4X) ``Advanced numerical methods for time dependent parametric partial differential equations with applications\rq\rq. We acknowledge financial support by Fondazione Cariplo and Fondazione CDP (Italy) under the project No.\ 2022-1895. This work has been written within the activities of GNFM and GNCS groups of INdAM (Italian National Institute of High Mathematics).

M.Z. acknowledges partial support of PRIN2022PNRR project No.P2022Z7ZAJ, European Union - NextGenerationEU and by ICSC - Centro Nazionale di Ricerca in High Performance Computing, Big Data and Quantum Computing, funded by European Union - NextGenerationEU.

\appendix
\section{Objective functions}
In Table \ref{tab:objective} we briefly recall the structure of the objective functions that have been considered in this paper. 

\begin{table}
\includegraphics[width=\linewidth]{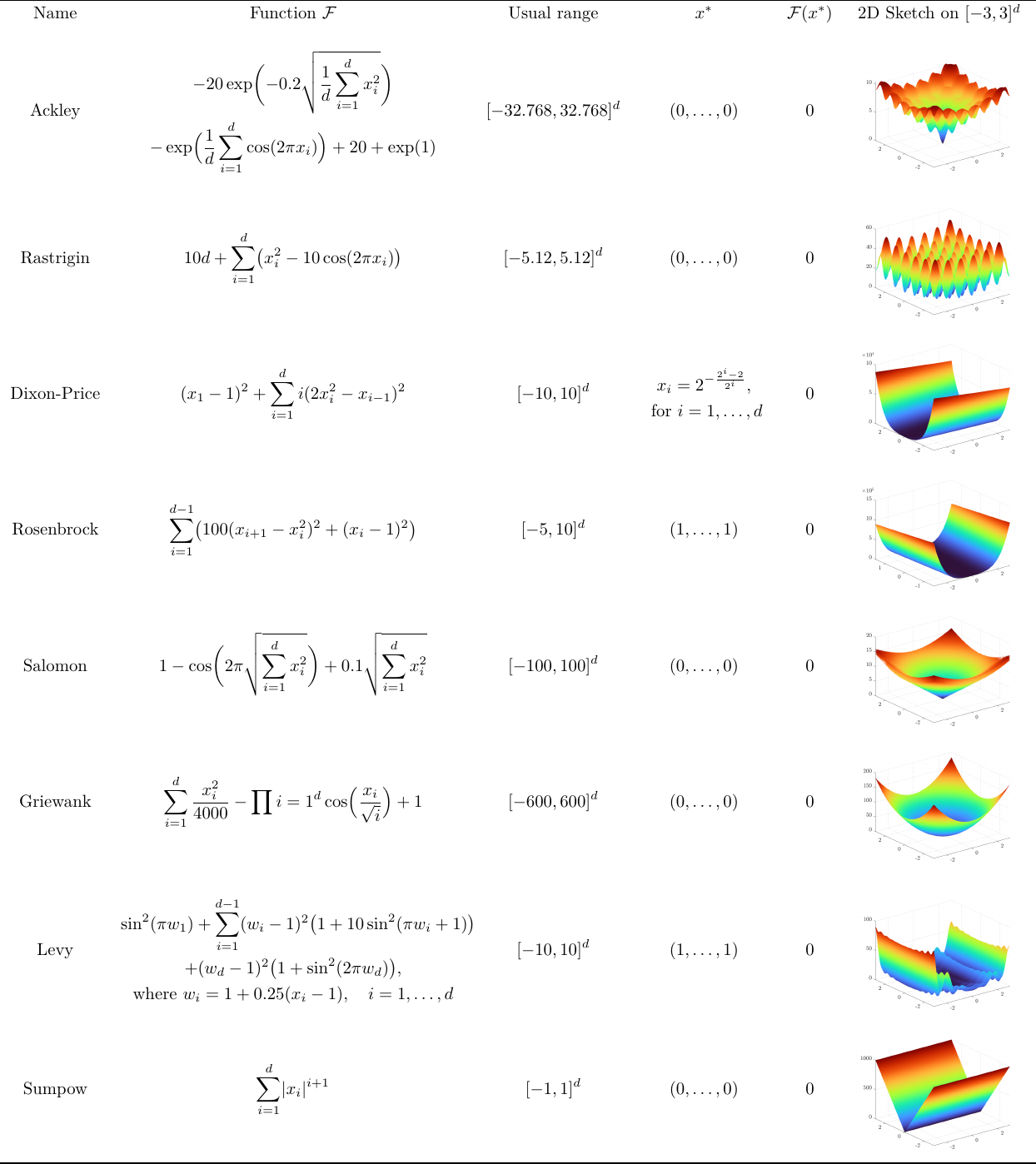}
\caption{Objective functions used for multidimensional optimization.}
\label{tab:objective}
\end{table}

\end{document}